\documentclass{article}
\usepackage[utf8]{inputenc}
\usepackage[top=1in,bottom=1.25in,left=1.25in,right=1.25in]{geometry}
\newcommand{\be}{\begin{equation}}
\newcommand{\ee}{\end{equation}}

\usepackage{scalerel,stackengine}
\stackMath
\newcommand\reallywidehat[1]{%
\savestack{\tmpbox}{\stretchto{%
  \scaleto{%
    \scalerel*[\widthof{\ensuremath{#1}}]{\kern-.6pt\bigwedge\kern-.6pt}%
    {\rule[-\textheight/2]{1ex}{\textheight}}
  }{\textheight}%
}{0.5ex}}%
\stackon[1pt]{#1}{\tmpbox}%
}
\usepackage{amsmath}
\usepackage{amsfonts}
\usepackage{amssymb}
\usepackage{graphicx}
\usepackage[dvipsnames]{xcolor}

\definecolor{LightGray}{rgb}{0.80,0.80,0.85}

\newcommand{\revTwo}[1]{{#1}}

\newcommand{\imunit}{\mathrm{i}}
\newcommand{\euler}{\mathrm{e}}
\newcommand{\Real}{\mathbb{R}}
\newcommand{\Complex}{\mathbb{C}}
\newcommand{\Integer}{\mathbb{Z}}
\newcommand{\Leb}{L}

\newcommand{\vect}[1]{\boldsymbol{#1}}
\newcommand{\fourier}[1]{\hat{#1}}
\newcommand{\bessel}[1]{\lowercase{#1}}
\newcommand{\spharm}[1]{\lowercase{#1}}
\newcommand{\wfourier}[1]{\widehat{#1}}
\newcommand{\rotation}{{\cal R}}
\newcommand{\translation}{{\cal T}}
\newcommand{\htranslation}{\widehat{\cal T}}
\newcommand{\sort}{{\text{\cal Sort}}}
\newcommand{\transpose}[1]{{#1}^{\intercal}}
\newcommand{\vx}{{\vect{x}}}
\newcommand{\vk}{{\vect{k}}}
\newcommand{\vd}{{\vect{\delta}}}
\newcommand{\vu}{{\vect{u}}}

\newcommand{\cX}{{\cal X}}
\newcommand{\cY}{\hat{\cal X}}
\newcommand{\cZ}{{\cal Z}}

\newcommand{\cN}{{\cal N}}
\newcommand{\cQ}{{\cal Q}}

\newcommand{\fA}{{\fourier{A}}}
\newcommand{\fB}{{\fourier{B}}}
\newcommand{\xF}{{        {F}}}
\newcommand{\fF}{{\fourier{F}}}
\newcommand{\xG}{{        {G}}}
\newcommand{\fG}{{\fourier{G}}}
\newcommand{\fsigma}{{\fourier{\sigma}}}

\newcommand{\fS}{{\fourier{S}}}
\newcommand{\fSlice}{{\fourier{Sl}}}
\newcommand{\bA}{{\bessel{A}}}

\newcommand{\bG}{{\bessel{G}}}

\newcommand{\sF}{{\spharm{F}}}
\newcommand{\sG}{{\spharm{G}}}

\newcommand{\vdmax}{{\vd_{\max}}}
\newcommand{\local}{{\text{\tiny local}}}
\newcommand{\vdlocal}{{\vd_{\local}}}
\newcommand{\kmax}{{K}}
\newcommand{\qmax}{{Q}}
\newcommand{\rmax}{{R}}

\newcommand{\tmax}{{T}}
\newcommand{\lmax}{{L}}
\newcommand{\mmax}{{M}}
\newcommand{\dx}{{\Delta x}}     
\newcommand{\dk}{{\Delta k}}     
\newcommand{\dpsi}{{\Delta \psi}}     
\newcommand{\Ccost}{C}
\newcommand{\Ckern}{\boldsymbol{C}}

\newcommand{\argmax}{\operatorname*{arg\,max}}
\newcommand{\bigO}{\mathcal{O}}

\newcommand{\nimage}{N_{A}}
\newcommand{\ntemplate}{N_{S}}
\newcommand{\nctf}{N_{\text{\tiny CTF}}}

\newcommand{\npixel}{N}
\newcommand{\nrank}{H}
\newcommand{\signal}{{\text{\tiny signal}}}
\newcommand{\noise}{{\text{\tiny noise}}}

\newcommand{\true}{{\text{\tiny true}}}
\newcommand{\estimated}{{\text{\tiny estim}}}
\newcommand{\empirical}{{\text{\tiny empir}}}
\newcommand{\initial}{{\text{\tiny init}}}
\newcommand{\random}{{\text{\tiny rand}}}
\renewcommand{\th}{{\text{\tiny th}}}
\newcommand{\opt}{{\text{\tiny opt}}}

\newcommand{\est}{{\text{\tiny est}}}

\newcommand{\tol}{{\text{\tiny tol}}}

\newcommand{\OracleAM}{{\text{\tiny OAM}}}
\newcommand{\Oracle}{{\text{\tiny Oracle}}}
\newcommand{\OracleAMStable}{{\text{\tiny OAM}}}
\newcommand{\RelionDeNovo}{{\text{\tiny RDN}}}
\newcommand{\BayesianInference}{{\text{\tiny BI}}}
\newcommand{\OracleBayesianInference}{{\text{\tiny OBI}}}
\newcommand{\OracleBayesianInferenceStable}{{\text{\tiny OBI}}}
\newcommand{\EMPM}{{\text{\tiny EMPM}}}
\newcommand{\masked}{{\text{\tiny mask}}}
\newtheorem{definition}{Definition}
\newtheorem{remark}{Remark}

\newcommand{\hk}{\hat{\boldsymbol{k}}}
\newcommand{\kpolara}{\theta}
\newcommand{\kazimub}{\phi}

\newcommand{\epolara}{\beta}
\newcommand{\eazimub}{\alpha}
\newcommand{\egammaz}{\gamma}

\title{Robust {\em ab initio} solution of the cryo-EM reconstruction problem at low resolution with small data sets}

\author{Aaditya V. Rangan and Leslie Greengard
\thanks{Center for Computational Mathematics, Flatiron Institute and Courant Institute, NYU}}

\date{July, 2023}

\begin{document}

\maketitle

\begin{abstract}
Single particle cryo-electron microscopy has become a critical tool in structural biology over the last decade, able to achieve atomic scale resolution
in three dimensional models from hundreds of thousands of (noisy) two-dimensional projection views of particles frozen at unknown orientations.
This is accomplished by using a suite of software tools to (i) identify particles in large micrographs, (ii) obtain low-resolution reconstructions, (iii) refine those low-resolution structures, and (iv) finally match the obtained electron scattering density to the constituent atoms that make up the macromolecule or macromolecular complex of interest.

Here, we focus on the second stage of the reconstruction pipeline:
obtaining a low resolution model from picked particle images.  Our goal is to create an algorithm that is capable of {\em ab initio} reconstruction from small data sets (on the order of a few thousand selected particles). 
{More precisely, we propose an algorithm that is robust, automatic, 
and fast enough that it can potentially be used to assist in the assessment of 
particle quality as the data is being generated during the microscopy experiment.}
\end{abstract}

\section{Introduction}

In single-particle cryo-electron microscopy (cryo-EM), a purified 
preparation of particles (proteins or macro-molecular complexes) 
is frozen in a thin sheet of ice, with each particle held in an unknown
orientation. Using a weak scattering approximation, the image
obtained can be viewed as a noisy projection of the particle's
electron scattering density, which we denote by $\xF(\vx)$,
with optical aberrations modeled by what 
is called a contrast transfer function (CTF)
\cite{Milne2013,Cheng2015}. The task at hand is to take the individual
particle images, determine their orientations, and recover a high 
resolution characterization of $\xF(\vx)$. This is typically accomplished
through the Fourier slice theorem 
\cite{Epstein,Natterer}, which states that
the two-dimensional Fourier transform of each
projection image is a slice through  $\fF(\vk)$,
the {\em three-dimensional} Fourier transform of $\xF(\vx)$
(modulated by the CTF). 

The specific slice obtained from a single image 
is determined by the orientation of the
individual particle. With many particles at many distinct orientations,
$\fF(\vk)$ is well-sampled and $\xF(\vx)$ can be obtained through the
inverse Fourier transform. Thus, if the orientations are known,
the reconstruction problem is linear and easily solved.
The orientations, however, are not known, and
the reconstruction problem is both nonlinear and non-convex.
In practice, not only are the orientations unknown, but the particles
must be correctly centered in order to invoke the Fourier slice theorem.
The unknown displacement vector must also be learned for 
each particle image.
We postpone a proper discussion of the mathematics to 
section \ref{sec:overview} and the appendices (Supporting Information).

From an experimental perspective, the initial stages of cryo-EM involve some 
hours of imaging, resulting in the collection of multiple micrographs,
each containing a number of particle images.
{After sufficiently many micrographs are collected, the next steps 
involve particle picking 
(identifying individual particles in the micrographs),
and the construction of a low- to medium-resolution molecule for further analysis
and refinement. The latter step typically involves
class averaging - finding particles that appear to have the same orientation
and taking the average in order to improve the signal-to-noise ratio, but some 
algorithms proceed from the raw particle images themselves (as we do here).}

{
Despite recent advances in sample preparation, imaging and analysis, the cost of data collection can still be 
quite significant, and the availability of microscope time is a limited resource.
In order for the picked-particle images to be useful, they must meet certain criteria: each should contain a single well-centered particle with a minimum of clutter, and the full set must cover a wide variety of viewing angles, sufficient to determine the full three-dimensional structure.
Assuming that these criteria are met, these images (hopefully the least noisy and cleanest projection views) 
can then be used in downstream analyses.}
If the sample or the images are not of sufficient quality, it may be necessary to terminate the current experiment and repeat with a new sample, hoping to produce higher-quality particle images the next time.

To reduce the time and cost of such experiments, there are many quality-control tasks that would be useful to carry out as early as possible.
These include:
\begin{enumerate}
\item determining whether any picked-particle image contains a reasonably well-centered and isolated molecule.
\item estimating the viewing angles of the images (and checking that the range of viewing angles is sufficient to reconstruct the molecule).
\item estimating the displacements/shifts required to more precisely center the images.
\item estimating the image quality and the number of particles that will be needed to reconstruct the desired molecule at the desired resolution.
\item estimating the molecular structure itself at low-resolution
(which feeds back to the estimates 1-4) 
\end{enumerate}
In current practice, many of these issues are addressed after carrying
out class averaging, which typically requires dozens of particles per class 
and tens of thousands of picked particles, {although a number of reconstruction algorithms work directly on the raw data.}

In this paper, we describe an algorithm which can serve as a complement to more standard strategies. 
Our algorithm operates on a relatively small set of picked-particle images (say, a few thousand), and attempts to directly reconstruct 
(without class averaging)
a low-resolution version of the imaged molecule, along with estimates of the viewing-angles, shifts and quality of the images involved.
As we demonstrate below, our method is significantly faster, more accurate and more reliable than {many} existing 
de-novo reconstruction algorithms.
We hope that it is robust enough to permit the assessment of data quality as the experiment is being carried out.
Moreover, our algorithm runs quickly and automatically with no `tuning' required, allowing for multiple small pools of approximately 1000 picked-particles to be independently assessed without supervision
(each corresponding to $\sim 30$ micrographs).
This would help set the stage for more robust statistical validation techniques such as 
bootstrap or jack-knife estimation.

{
\begin{remark}
We should emphasize that the
starting point of our experiments here is a collection of picked particle images 
(obtained from the EMPIAR database), not
the raw movies or entire micrographs that are sometimes deposited as well.
Particle picking is itself a non-trivial task.
Restricting our attention to molecules for which particles are curated
allows us to focus exclusively on the problem of reconstruction from a small
data set. 
\end{remark}}

At its core, our algorithm is a modification of the simplest expectation-maximization method.
We use the framework of classical iterative refinement - alternating between a reconstruction step and an alignment step.
In the reconstruction step, we use the current estimate of the viewing angles for each particle/image and their displacements to approximate the molecule. 
In the alignment step, we use the current estimate of the molecule to better fit the viewing angles and displacements for each image. We will refer to this approach as {\em alternating minimization} (AM).

It is well-known that standard AM does poorly in this setting; an initialization with randomized viewing angles typically fails to converge to a useful low-resolution approximation of the molecule.
To address these issues and to guide AM to a better structure, we introduce two changes: 
(i) we compress and denoise the individual images by computing their 
two-dimensional Fourier transform on a polar grid and focusing on carefully selected rings in the transform domain with maximal information content, and (ii) we increase the entropy of the viewing-angle distribution used during the reconstruction step (in a manner to be described in detail).
These modifications both reduce the computational cost and improve the robustness of AM in the low-resolution low-image-number regime.

In the numerical experiments below,
we show that our modified method ({\em alternating minimization with 
entropy maximization and principal modes} or `EMPM') 
performs surprisingly well on several published data sets. 
Our reconstructions are often on par with the best results one could expect from an `oracle-informed' AM (i.e., standard AM starting with the `ground-truth' viewing-angles and displacements for each image).
Moreover, our results are similar to (or better than) what one might expect from a reconstruction using `optimized' Bayesian inference (i.e., using a ground-truth informed estimated noise level).
Importantly, our strategy is efficient, running about twenty times 
faster than standard Bayesian inference or the 
de-novo reconstruction in the widely used, open source package Relion 
\cite{Kimanius2016,Kimaniusetal2021,Scheres2012}.
{We also make direct comparisons with two other open source 
methods: the neural network approach in
cryoDRGN \cite{ZBBD21} and the recent low-resolution reconstruction pipeline in Aspire \cite{AspireWebsite}.
We have not made a direct comparison with the stochastic gradient descent approach in 
the GPU-based cryoSPARC package \cite{Punjani2017}.} 

{
In the remainder of this paper,
we briefly summarize the EMPM method and demonstrate its performance on a 
variety of data sets.
The details of our method, as well as some additional motivating examples, 
are deferred to the Supporting Information.
It is currently implemented for multi-core CPUs, but we believe it should be 
straightforward to implement on GPUs as well, and to integrate into any of the 
existing software packages for cryo-EM.}

\section{Overview of the EMPM iteration} \label{sec:overview}

Our general goal is to solve the following problem: 
given a set of $\nimage$ picked particle images (2D projections) 
from a handful of 
micrographs (each with their own CTF), 
estimate the unknown viewing-angles $\tau_{j}\in SO(3)$ and displacements $\delta_{j}\in \Real^{2}$ for each image, 
as well as the underlying volume $\fF(\vk)$,
where $\vk\in\Real^{3}$ denotes a point in the Fourier transform domain.
In addition to determining the volume $\fF(\vk)$, we would also like 
an estimate of the viewing-angle distribution $\mu_{\tau}$ and 
displacement distribution $\mu_{\vd}$ (taken across the images).
If the estimates of the volume $\fF(\vk)$ and the distributions 
$\mu_{\tau}$, $\mu_{\vd}$ are 
accurate estimates of the ground truth, they
will clearly be useful for quality control, determining (for example) 
whether the experimental conditions are suitable for larger scale data
collection and subsequent high resolution reconstruction.

The main tool we employ is a form of {\em alternating minimization} (AM), 
which carries out some number of iterations of the following two 
step procedure:
\begin{description}
\item[Reconstruction:] Given the estimated viewing-angles $\{\tau_{j}^{\estimated}\}$ and displacements $\{\vd_{j}^{\estimated}\}$ for $\nimage$ 2d picked-particle images $\{\fA_{j}\}$, reconstruct the 3D molecular volume $\fF^{\estimated}$ (using the Fourier slice theorem and least squares inversion of the Fourier transform).
\item[Alignment:] Given the reconstruction of the 3D molecule 
$\fF^{\estimated}$, update the estimates $\{\tau_{j}^{\estimated}\}$ 
and $\{\vd_{j}^{\estimated}\}$.
\end{description}

In a standard AM scheme, the reconstruction step solves 
for $\fF^{\estimated}$ using a least-squares procedure, and the 
alignment step assigns the $\{\tau_{j}^{\estimated}\}$ and 
$\{\vd_{j}^{\estimated}\}$ using some kind of maximum-likelihood 
procedure (see sections \ref{sec_Volume_Reconstruction} and \ref{sec_Image_Alignment} in the Supporting Information).
While there are many more sophisticated strategies for molecular 
reconstruction (see, for example, \cite{Punjani2017,ZBBD21}), 
we focus on AM for 
a few reasons:
first and foremost is its simplicity, transparency and ease of automation.
AM has modular steps, few parameters, and is easily
interpretable.

At the same time, however, standard AM does have serious limitations.
First, it can be computationally expensive, particularly if a broad range 
of frequencies and/or translations are considered in the alignment step. 
Second, AM can fail to converge to the correct solution, especially
during the initial iterations, when the initial guess is far from the 
true molecule. 

To address these two issues, we propose a modified version of AM 
which we refer to as {\em alternating minimization with entropy-maximization
and principal modes} (EMPM). The two main modifications we introduce are: 
\begin{description}
\item[Maximum entropy alignment:] We use both the traditional 
maximum likelihood alignment strategy, which finds the best estimate 
for each image from the current molecular reconstruction, and its
{\em opposite}. That is, we also use a
`maximum-entropy' alignment strategy: for each uniformly chosen
viewing angle, we
find the experimental image which best matches the corresponding projection
from the current structure. The subsequent reconstruction step makes use of
this uniform distribution of viewing angles 
{\em whether or not it corresponds to the true distribution}.
Use of this counter-intuitive step is important in stabilizing 
our low resolution reconstruction from a small set of images
(see section \ref{sec_Image_Alignment} in the Supporting Information).
\item[Principal-mode projection:] 
Each image is transformed onto a polar grid in Fourier space.
The data on each circle in 
the transform domain at a fixed modulus $|\vk|$ will be referred to as a
{\em Fourier ring}.
Rather than considering all Fourier rings when comparing
2D projections and image data,
we consider only a handful of carefully-chosen combinations of such rings,
which we refer to as {\em principal modes}.
This yields significant compression of the data 
(and computational efficiency) while retaining information useful for 
alignment (see section \ref{sec_Principal_image_rings} in the 
Supporting Information). 
\end{description}

In the EMPM pipeline, we first estimate the principal-modes using 
the images themselves. 
Using these `empirical' principal modes, we proceed with AM,
alternating between maximum-likelihood and maximum-entropy alignment.
Once a preliminary volume estimate has been obtained, 
we use that estimate to recalculate the principal modes.
We then restart AM using the new principal images and, again,
alternate between maximum-likelihood and maximum-entropy alignment 
until we have a final volume estimate.
The algorithm is summarized below (see Supporting Information for a more detailed description):
\begin{enumerate}
\item Use the experimental images $\fA_{j}$ to calculate 
`empirical' principal modes $U^{\empirical}$.
\item Compute the {\em principal images}, defined to be
the projection of the images onto principal modes, 
$P(U^{\empirical},\fA_{j})$.
\item Set random initial viewing angles $\tau_{j}^{\empirical}$ and set 
initial displacements $\vd_{j}^{\empirical}$ to $0$.
\item Carry out 32 iterations of EMPM on the principal images. 
\begin{enumerate}
\item Reconstruct $\fF^{\empirical}$ using $\tau_{j}^{\empirical}$ and $\vd_{j}^{\empirical}$.
\item Use maximum-likelihood alignment to update $\tau_{j}^{\empirical}$ and $\vd_{j}^{\empirical}$.
\item Reconstruct $\fF^{\empirical}$ using $\tau_{j}^{\empirical}$ and $\vd_{j}^{\empirical}$.
\item Use maximum-entropy alignment to update $\tau_{j}^{\empirical}$ and $\vd_{j}^{\empirical}$.
\end{enumerate}
\item Use $\fF^{\empirical}$ to calculate `estimated' principal modes $U^{\estimated}$.
\item Set initial viewing angles $\tau_{j}^{\estimated}=\tau_{j}^{\empirical}$ and displacements $\vd_{j}^{\estimated}=\vd_{j}^{\empirical}$.
\item Carry out 32 iterations of EMPM on the new principal images 
$P(U^{\estimated},\fA_{j})$.
\begin{enumerate}
\item Reconstruct $\fF^{\estimated}$ using $\tau_{j}^{\estimated}$ and $\vd_{j}^{\estimated}$.
\item Use maximum-likelihood alignment to update $\tau_{j}^{\estimated}$ and $\vd_{j}^{\estimated}$.
\item Reconstruct $\fF^{\estimated}$ using $\tau_{j}^{\estimated}$ and $\vd_{j}^{\estimated}$.
\item Use maximum-entropy alignment to update $\tau_{j}^{\estimated}$ and $\vd_{j}^{\estimated}$.
\end{enumerate}
\end{enumerate}
Note that, each time we update the displacements $\vd_{j}^{\empirical}$ and $\vd_{j}^{\estimated}$, we bound the maximum permitted displacement magnitude by $\vdmax$ (typically $\vdmax\sim 0.1$ in our dimensionless units).
In our studies with small amounts of moderately noisy 
experimental data, we have observed that
the EMPM strategy is more effective than 
Bayesian inference or maximum-likelihood alignment
alone - even with probabilistic alignment where a distribution of viewing angles is assigned to each image.
Roughly speaking, we believe this is due the fact that the incorrect viewing angles assigned in the 
maximum-entropy step are washed out by destructive interference, while the
correctly assigned ones reinforce the true structure. 
One simple model problem which illustrates these points is discussed in the Supporting Information (see sections \ref{sec_MRA} and \ref{sec_MSA}).
It is worth repeating, however, that we only propose this approach for obtaining a low-resolution initial guess. Existing refinement methods are still crucial for high resolution.

Once we have run our EMPM algorithm, we have the estimated volume $\fF^{\estimated}(\vk)$ as output.
Using $\fF^{\estimated}(\vk)$ as an approximation of $\fF(\vk)$, we can readily compute other quantities that may be of interest, such as 
estimates of the true viewing angle distribution in the data set and other types of correlations.
We will use these quantities in our numerical experiments below to compare the results of EMPM to other methods, such as Relion's de-novo molecular reconstruction and full Bayesian inference.

\subsection{Templates, correlations, and ranks} \label{sec:correlations}

For each centered molecular reconstruction $\fF$, and set of 
Euler angles $\tau = (\alpha,\beta,\gamma)$, we define 
the inner product $\cX$ and normalized inner product $\cY$ 
of the reconstruction with the ``ground-truth" model
$\fF^{\true}$ by 
\begin{eqnarray}
\cX(\tau;\fF,\fF^{\true}) & = & \langle \rotation(\tau)\circ \fF , \fF^{\true} \rangle \\
\cY(\tau;\fF,\fF^{\true}) & = & \frac{\cX(\tau;\fF,\fF^{\true})}{\|\fF\|\cdot\|\fF^{\true}\|} \, . \label{Eq_innerprod}
\end{eqnarray}
Here $(\alpha,\beta)$ denotes the viewing angle of an image with
$\alpha$ the azimuthal angle and $\beta$ the polar angle. $\gamma$ is
the final in-plane rotation (see section \ref{sec_Volume_notation}).
$\rotation(\tau)\circ \fF$ is the volume obtained from 
$\fF$ using the rotation (the element of $SO(3)$) defined by $\tau$. 
We define the volumetric correlation $\cZ(\fF,\fF^{\true})$ by
\begin{equation}
\cZ(\fF,\fF^{\true}) = \max_{\tau} \cY(\tau;\fF,\fF^{\true}) \text{.} 
\label{Eq_volume_correlation}
\end{equation}
This measurement $\cZ(\fF,\fF^{\true})$ represents the optimal 
correlation (over all alignments) between the two 
volumes $\fF$ and $\fF^{\true}$.

In order to compare our results with Relion's de-novo reconstruction, 
we will also need to calculate the `masked' correlation 
$\cZ^{\masked}(\xF,\xF^{\true})$, defined as:
\begin{eqnarray}
\cX^{\masked}(\tau;\xF,\xF^{\true}) & = & \langle \rotation(\tau)\circ \xF\cdot\xG , \xF^{\true} \rangle \\
\cY^{\masked}(\tau;\xF,\xF^{\true}) & = & \frac{\cX^{\masked}(\tau;\xF,\xF^{\true})}{\|\xF\cdot\xG\|\cdot\|\xF^{\true}\|} \\
\cZ^{\masked}(\xF,\xF^{\true}) & = & \max_{\tau} \cY^{\masked}(\tau;\xF,\xF^{\true}) \text{,} 
\label{Eq_volume_correlation_masked}
\end{eqnarray}
where $\xG$ is a sharp spherical mask (i.e., the indicator-function of a ball centered at the origin) with a radius determined by the support of $\xF^{\true}$.

We will also need to define inner products and correlations of single
images and volumetric slices. For this, given a Fourier space volume
$\fF$, we generate noiseless 2D projections (in the Fourier domain) 
that incorporate the CTF, which we refer to as {\em templates}.
We denote these by 
$\fS(\tau , \vd ; CTF ; \fF)$,
where $\tau$ is the specific set of Euler angles 
and $\vd$ is the displacement.
When assuming the displacement is zero, we will write these templates
as $\fS(\tau; CTF ; \fF)$,
Formally, we have
\begin{eqnarray}
\fS(\tau ; CTF ; \fF) & := & CTF(\vk)\odot\left[\fSlice\circ\rotation(\tau)\circ \fF\right](\vk) \text{.}
\label{eq_translated_CTF_modulated_template} 
\end{eqnarray}
where $CTF(\vk)$ is the given pointwise value of the CTF in Fourier space
and $\fSlice$ denotes taking the equatorial slice in the Fourier domain.
These CTF-modulated templates are used to calculate the inner products:
\begin{eqnarray}
\cX(\tau , \vd ; \fA_{j} ; CTF_{j} ; \fF ) 
  & = & \langle \ \fS(\tau ; CTF_{j} ; \fF) \ , \ \translation(-\vd,\vk)\odot{ \fA_{j}} \ \rangle \text{,}
\end{eqnarray}
where $\translation(+\vd,\vk) = e^{i \vd \cdot \vk}$ is the Fourier signature of displacement  by 
the vector $\vd$.
These are used, in turn, to calculate the correlations:
\[
\cY(\tau , \vd ; \fA_{j} ; CTF_{j} ; \fF) = \frac{\cX(\tau , \vd ; \fA_{j} ; CTF_{j} ; \fF)}{ \| \fS(\tau ; CTF_{j} ; \fF) \| \cdot \| \fA_{j} \| } \text{.}
\]
For convenience, we generally reduce this to a function of the 
viewing direction 
$\hk = \left(\eazimub,\epolara\right)$, the first two components of 
$\tau$, alone:
\begin{equation}
\cZ(\eazimub , \epolara ; \fA_{j} ; CTF_{j} ; \fF) = \max_{\egammaz,\vd} \cY(\tau,\vd ; \fA_{j} ; CTF_{j} ; \fF) \text{.}
\label{eq_cZ}
\end{equation}
The parameters $\egammaz$ and $\vd$ are associated with a particular choice of $\eazimub$ and $\epolara$ implicitly through the $\argmax$.
When the context is clear, we will often abuse notation and refer to 
the correlations in \eqref{eq_cZ} by
$\cZ(\eazimub , \epolara ; \fA_{j})$.

Finally, for each image $\fA_j$ and a set of $\ntemplate$ templates
defined by the viewing directions $(\alpha_l,\beta_l), \ l=1,\dots,\ntemplate$, 
we can order the values 
\begin{equation}
\cZ(\eazimub_l , \epolara_l ; \fA_{j} ; CTF_{j} ; \fF), \ l = 1,\dots,\ntemplate.
\label{templaterankdef}
\end{equation}
This yields what we will refer to as the {\em template rank} for the image.

Conversely, for each template defined by $(\eazimub_l , \epolara_l)$, 
one can order the values
\begin{equation}
\cZ(\eazimub_l , \epolara_l ; \fA_{j} ; CTF_{j} ; \fF), \ j = 1,\ldots,\nimage.
\label{imagerankdef}
\end{equation}
This defines the  {\em image rank} for that template.
Using the terminology of linear algebra, we can define the 
$\ntemplate\times\nimage$ matrix 
$\cZ$ with
\[ \cZ_{lj} = \cZ(\eazimub_{l},\epolara_{l};\fA_{j}).
\] 
Template-ranking then corresponds to sorting the rows of $\cZ$ while
image-ranking corresponds to sorting the columns.

\section{Results}

We will consider a number of examples in this paper, drawn from the 
electron-microscopy image archive (EMPIAR) - namely, {\em trpv1},
the capsaicin receptor (the 
transient receptor potential cation channel subfamily V member 1) 
[EMPIAR-10005, emd-5778],
{\em rib80s}, the Plasmodium falciparum 80S ribosome bound to emetine
[EMPIAR-10028, emd-2660],
{\em ISWINCP}, the ISWI-nucleosome complex
[EMPIAR-10482, emd-9718],
{\em ps1}, the pre-catalytic spliceosome 
[EMPIAR-10180],
{\em MlaFEDB}, a bacterial phospholipid transporter
[EMPIAR-10536, emd-22116],
{\em LetB1}, the lipophilic envelope-spanning tunnel B protein from E. coli
[EMPIAR-10350, emd-20993],
{\em TMEM16F}, a Calcium-activated lipid scramblase and ion channel in 
digitonin with calcium bound
[EMPIAR-10278, emd-20244], and
{\em LSUbl16dep}, 
L17-Depleted 50S Ribosomal Intermediates, using an average of the published major classes as a reference
[EMPIAR-10076].
To limit the number of picked particles,
we carry out our numerical experiments on 
subsets of $\nimage=1024$ images, spanning a small number 
$\nctf$ of micrographs.

We set $\kmax=48$ as the bandlimit for our low-resolution reconstruction, corresponding to a `best' possible resolution of about $20$\AA. 
With these image pools, we construct a 
variety of different `reference volumes' for $\fF(\vk)$.
\begin{description}
\item[Ground truth:] First, we center and restrict the published volume to the band-limit $\kmax$. 
This defines our `ground truth' volume $\fF^{\true}(\vk)$.
\item[Oracle:] Next, we measure the correlations 
$\cZ(\tau;\fA_{j};CTF_{j};\fF^{\true})$ for a sufficiently sampled
set of Euler angles $\tau$ and displacements and use maximum-likelihood to align the images to $\fF^{\true}(\vk)$. This assigns `ground-truth' 
viewing angles $\tau_{j}^{\true}$ and displacements $\vd_{j}^{\true}$ to 
each image. More precisely, we start by calculating $\cZ$ using a small 
displacement disc of $\vdlocal=0.03$, and then iterate multiple times, 
with $\vdmax=0.25$ until the correlation between the reconstructed and true 
volumes plateaus (this usually takes less than $8$ iterations). We then 
take the final viewing angles and displacements to define $\tau_{j}^{\true}$ and $\vd_{j}^{\true}$. Given these ground-truth-based viewing angles and displacements, we use a least-squares reconstruction to compute $\fF^{\Oracle}(\vk)$.
\item[Oracle-AM (OAM):] Starting with the ground-truth 
viewing angles and displacements, as well as 
$\fF^{\Oracle}(\vk)$, we apply standard AM,
iterating until the resulting volume converges (this, again, usually
requires less than $8$ iterations). 
We define this final volume to be $\fF^{\OracleAM}(\vk)$. 
We expect the volume $\fF^{\OracleAM}$ to be close to the best 
reconstruction one could realistically achieve with standard AM. 
\item[Oracle-Bayesian (OBI):] 
Next, we apply Bayesian inference reconstruction 
(as described in \cite{Scheres2012,Scheres2012b}).
This reconstruction depends on (i) the initial starting volume, 
(ii) the estimated noise level $\sigma_{\BayesianInference}$, 
(iii) the resolution $\kmax$ and range of displacements $\vdmax$ considered,
and (iv) the error tolerance used when calculating the integrals required 
for estimating the posterior described in \cite{Scheres2012}. 
We fix $\kmax = 48$ and $\vdmax=0.10$ and set the error tolerance to the 
same value used in the rest of our numerical experiments 
(typically $\epsilon_{\tol}=0.01$). 
To create a reference volume associated with the ground truth, 
we use $\fF^{\true}$ as the initial starting volume.
The resulting volume -- denoted by 
$\fF^{\OracleBayesianInference}(\vk;\sigma_{\BayesianInference})$ -- will depend on the noise-estimate $\sigma_{\BayesianInference}$.
Note that, when $\sigma_{\BayesianInference}\rightarrow 0$ we will have $\fF^{\OracleBayesianInference}(\vk;\sigma_{\BayesianInference})\rightarrow\fF^{\OracleAM}(\vk)$, and when $\sigma_{\BayesianInference}\rightarrow \infty$ we will have $\fF^{\OracleBayesianInference}(\vk;\sigma_{\BayesianInference})$ converge to a spherically-symmetric distribution. 
In our experience the quality of the reconstructed volume depends quite sensitively on the choice of $\sigma_{\BayesianInference}$.
To build a good noise-estimate, we scan over several values of $\log(\sigma_{\BayesianInference})\in[-5,\ldots,0]$, choosing the $\sigma_{\BayesianInference}$ for which $\fF^{\OracleBayesianInference}(\vk;\sigma_{\BayesianInference})$ has the highest correlation with $\fF^{\true}$ after the first iteration.
Fixing $\sigma_{\BayesianInference}^{\opt}$ to be this optimal value, we define $\fF^{\OracleBayesianInference}(\vk):=\fF^{\OracleBayesianInference}(\vk,\sigma_{\BayesianInference}^{\opt})$ after iterating to convergence.
Thus, $\fF^{\OracleBayesianInference}$ is close to the best reconstruction possible using Bayesian inference, assuming that the optimal estimated noise-level $\sigma_{\BayesianInference}$ is known ahead of time.
\end{description}

In addition to the above reference volumes, we also create several
de-novo molecular reconstructions. RDN and EMPM below are 
oracle-free, while for Bayesian inference we use the optimal
noise level estimate $\sigma_{\BayesianInference}^{\opt}$ from
above.

\begin{figure}
    \centering
    \includegraphics[width=6.0in]{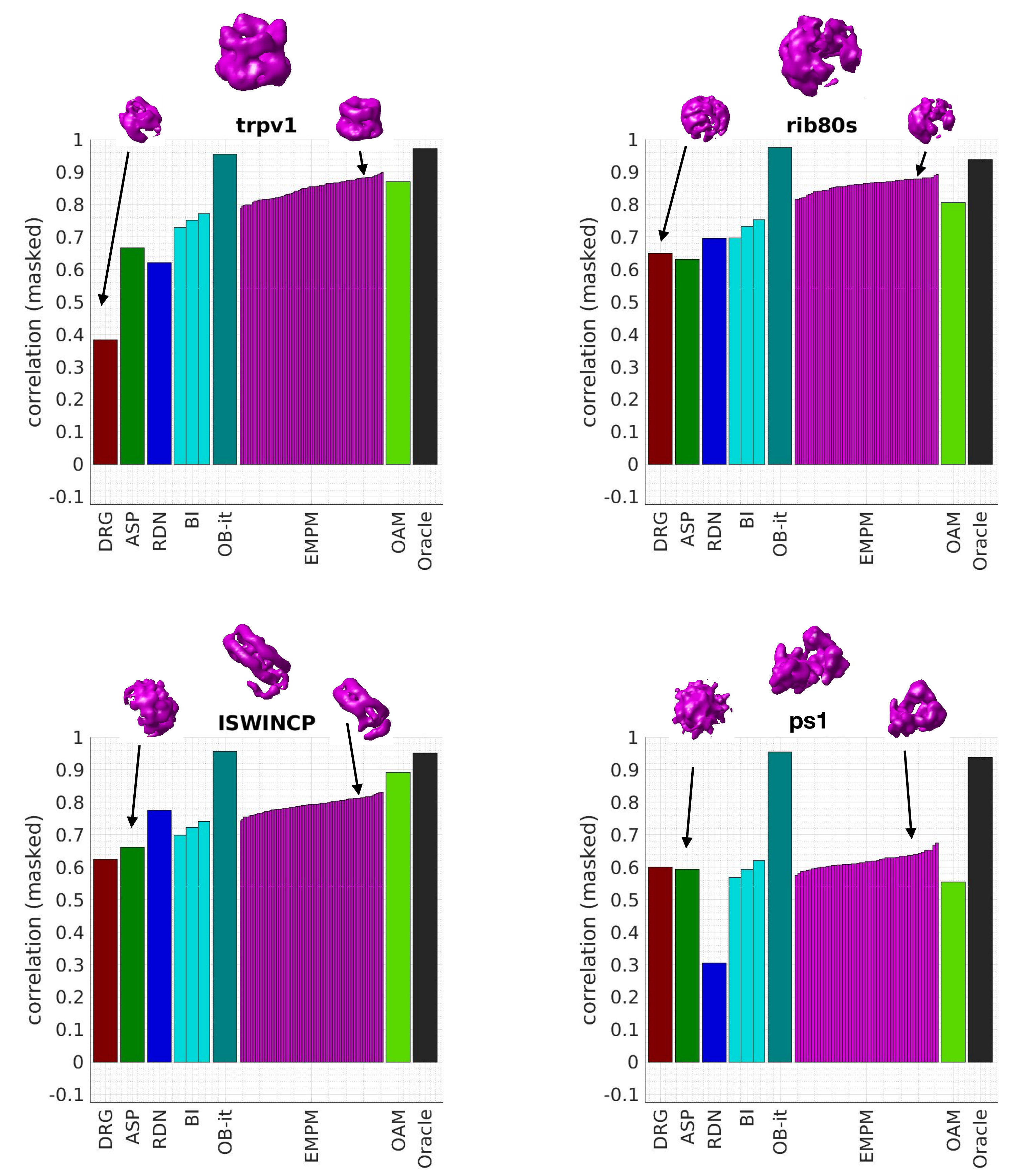} \\
    \caption{
{
    In this figure we illustrate the aligned correlation $\cZ^{\masked}$ between our various reconstructions and the ground-truth for 
{\em trpv1, rib80s, ISWINCP}, and {\em ps1} (see text).
    The different horizontal bars indicate the values of $\cZ^{\masked}$ for different reconstruction strategies (coded by color) with different trials from the same strategy adjacent to one another (sorted by $\cZ^{\masked}$).
    }
{At the top of each panel (centered) is the ``Oracle" low resolution
reconstruction for the indicated molecule. To the right in each panel, we show an EMPM 
reconstruction and to the left, we show a reconstruction obtained using one of the other pipelines
which achieves poorer correlation scores. Note that the $\cZ^{\masked}$ values
do tend to correlate with the image reconstruction quality.
{\em trpv1} has a highly non-uniform distribution of ground-truth Euler angles. The
distributions in the other cases are more uniform. (See Fig. \ref{Fig_view_angle}.)
}
    \label{Fig_collect_SET_ONE}}
\end{figure}

\begin{description}

\item[Bayesian Inference (BI):] We run the Bayesian inference reconstruction described in \cite{Scheres2012,Scheres2012b} to produce $\fF^{\BayesianInference}(\vk)$.
As mentioned above, this reconstruction depends on (i) the initial starting volume, (ii) the estimated noise-level $\sigma_{\BayesianInference}$, (iii) the resolution $\kmax$ and the range of displacements  $\vdmax$ considered, (iv) the error-tolerance used when computing the integrals in the posterior, and (v) the number of iterations applied. 
We fix the noise level $\sigma_{\BayesianInference}^{\opt}$ to be the value obtained in the Oracle-Bayesian method above.
We fix the resolution at $\kmax=48$, set $\vdmax=0.10$, set the error-tolerance to the same $\epsilon_{\tol}$ used in the other methods, and set the number of iterations at $16$. 
For the initial starting volume we assign each of the available images a viewing angle randomly and uniformly, and then use a least-squares 
reconstruction to generate a `random' initial starting volume. 
While random, this starting volume has roughly the same diameter and bandwidth as the true molecule.
We use the same procedure when initializing our strategy below. 
{
\item[DRG:] We run the open-source neural-network-based algorithm
cryoDRGN described in \cite{ZBBD21}, with the suggested default parameters for {\em ab initio}
reconstruction.
\item[ASP:] We run the open-source software package ASPIRE \cite{AspireWebsite},
which uses a common-lines-based algorithm, with the suggested default parameters
for {\em ab initio} reconstruction.
}
\item[RDN:] We run Relion's de-novo molecular reconstruction algorithm
(RDN) - a variant of Bayesian inference -
to produce $\fF^{\RelionDeNovo}$.
{(For speed, RDN uses only subsets of
images in its internal iterations.)}
We use the parameters listed in section \ref{sec_RelionDeNovo} in the Supporting Information.
We set the parameter {\tt \small --particle\_diameter} to be 
commensurate with the true molecular diameter, 
rounded up to the nearest $50$\AA. 
\item[EMPM:] We carry out the reconstruction using the strategy described in
the Methods section above, and discussed in more detail in 
the Supporting Information to produce $\fF^{\EMPM}(\vk)$. 
In addition to setting the maximum
frequency $\kmax$ and displacement $\vdlocal$, many of the computed 
quantities require a numerical tolerance $\epsilon_{\tol}$
(see section \ref{sec:notation}).
We fix $\kmax=48$, and scan over $\epsilon_{\tol}\in[0.001,0.010]$ and $\vdlocal\in[0.01,0.10]$. 
Our strategy begins with the same random initial viewing-angles for 
each image as used for Bayesian inference. 
As we shall see below (see, e.g., Fig. \ref{Fig_collect_SET_ONE}), our 
results are quite robust, provided that $\epsilon_{\tol}\leq 0.01$ 
and $\vdlocal\geq 0.01$ or so.
\end{description}

For {\em trpv1}, as well as {\em rib80s} and {\em ISWINCP}, the data sets are quite clean. 
The first $N_A = 1024$ picked-particle images
are well-centered, devoid of clutter and of high quality. 
Consequently, many of the reconstructions we attempted were successful 
and RDN reconstruction was also able to achieve 
$\cZ^{\masked} > 0.5$.
For {\em ps1, LSUbl16dep, MlaFEDB, LetB1}, and {\em TMEM16F}, however,
the individual images are not quite as clean. In the case of {\em ps1} and {\em TMEM16F},
many of the images are particularly noisy. 
For {\em LSubl16dep} (and to a lesser extent {\em ps1}), a notable fraction of the images are poorly centered. 
For {\em MlaFEDB} (and to a lesser extent {\em TMEM16F}), many of the images appear to contain fragments of other particles. For {\em LetB1}, the viewing-angle distribution for the images is highly concentrated around the equator. 
Our reconstructions for these molecules are not as successful. 
Nevertheless, we still recover volumes that are comparable in quality to Bayesian inference with an optimal estimated noise-level, and usually better than RDN reconstruction (all at a fraction of the computational time).

\begin{figure}
    \centering
    \includegraphics[width=6.0in]{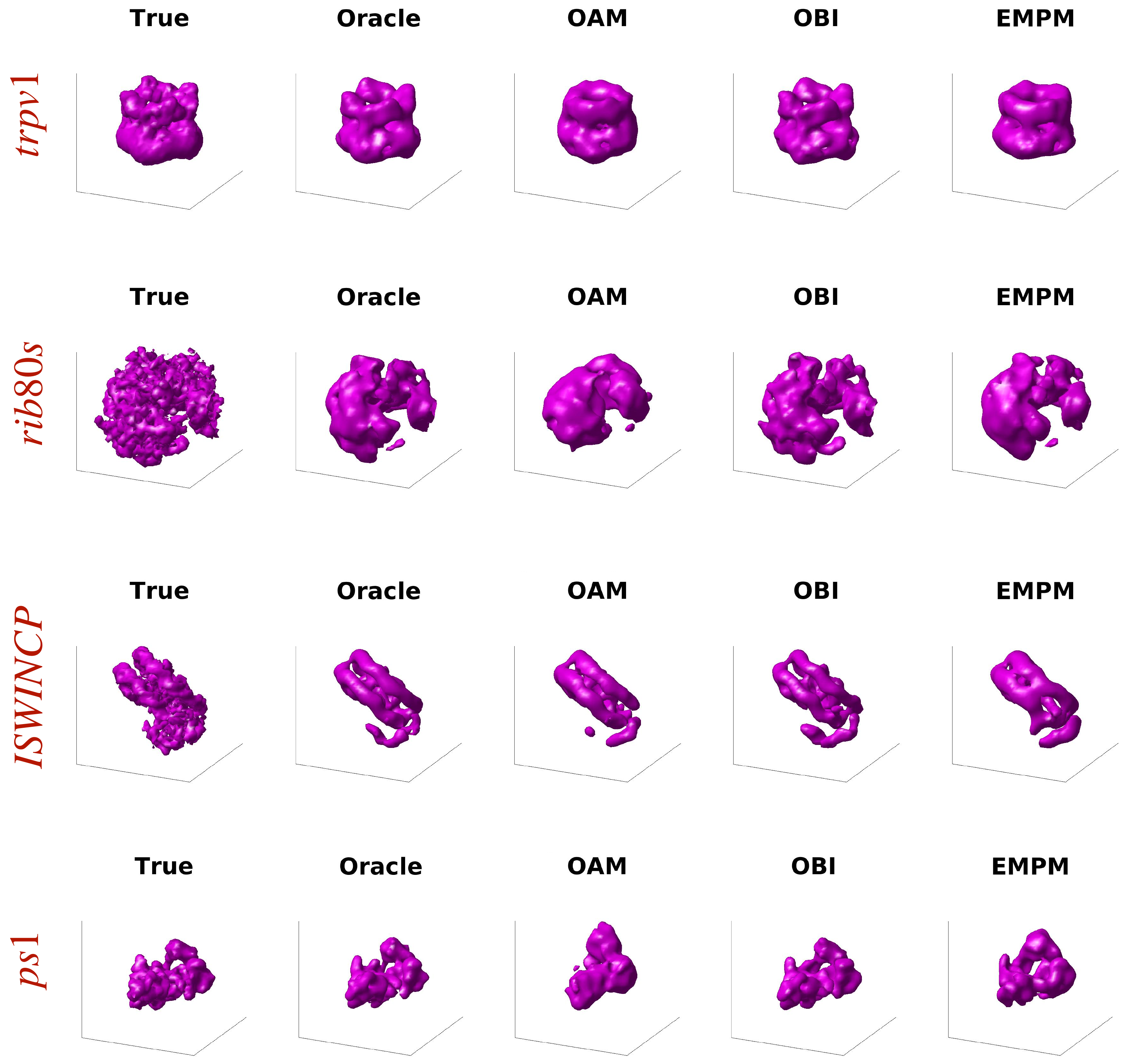} \\ 
    \caption{
{
    In this figure, we plot level sets from the volumes $\fF^{\true}$, $\fF^{\Oracle}$, $\fF^{\OracleAMStable}$, and $\fF^{\OracleBayesianInferenceStable}$, as well as one of the $\fF^{\EMPM}$ reconstructions (left to right, respectively).
    From the top down, the panels correspond to {\em trpv1},
    {\em rib80s}, {\em ISWINCP} and {\em ps1}.}
    \label{Fig_vol_SET_ONE}}
\end{figure}

\begin{figure}
    \centering
    \includegraphics[width=6.0in]{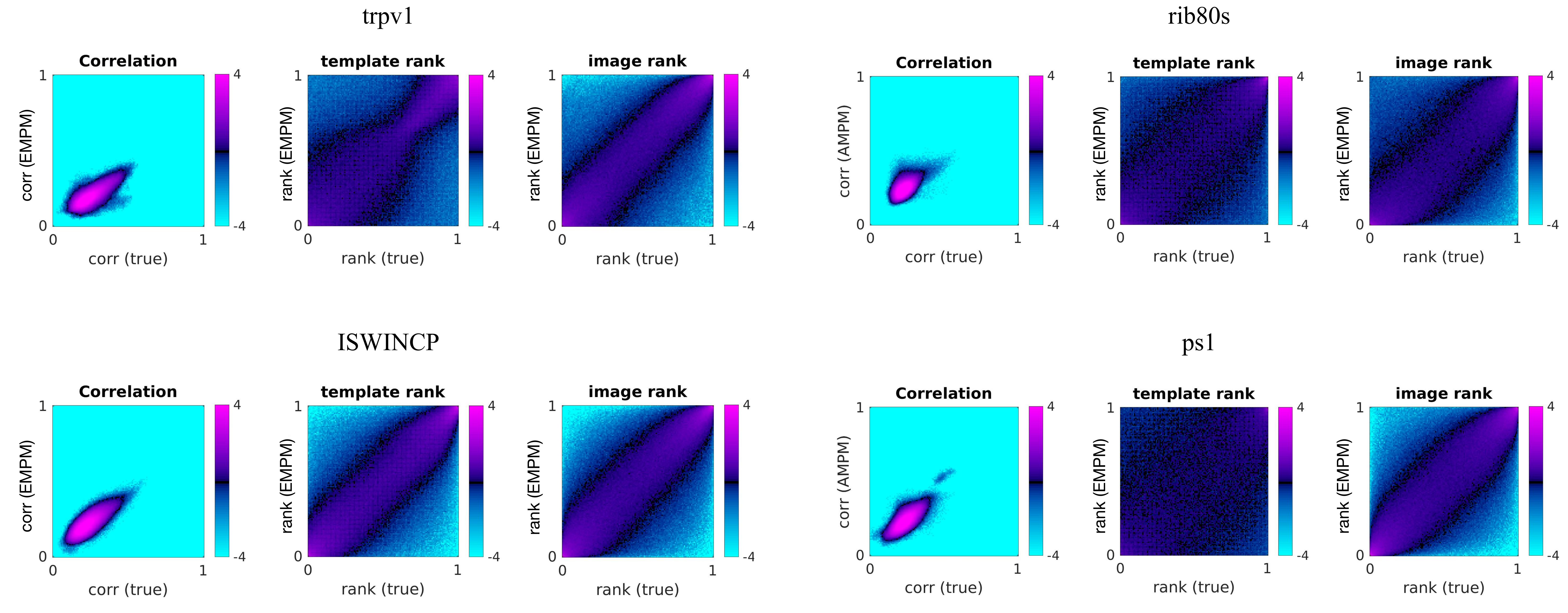} \\ 
    \caption{
    In this figure we illustrate the correspondence between estimated and ground-truth viewing angles for the reconstructions shown in Fig \ref{Fig_vol_SET_ONE}.
Three panels are shown for each molecule of interest. 
The left-most subplot
is a scatterplot of the image-template correlations. 
The middle subplot is a scatterplot of the template
ranks aggregated over all images, and
the right-most 
subplot is a scatterplot of the image ranks aggregated over all templates.
All are presented as heatmaps, with color corresponding to the
log-density of the joint-distribution.
(See text for further explanation.)
}
    \label{Fig_viewing_angle_distribution_SET_ONE}
\end{figure}

Results for four molecules are shown in Fig. \ref{Fig_vol_SET_ONE}.
Focusing on {\em trpv1} for the moment,
we see that $\fF^{\EMPM}$ captures the coarse features of $\fF^{\true}$, 
and is not too different from $\fF^{\OracleAMStable}$ in terms of overall quality.  Even though $\fF^{\EMPM}$ is a coarse estimate of $\fF^{\true}$, the reconstruction process has recovered a significant amount of information regarding the viewing-angles of each image.
To illustrate this point, we compare the image-template correlations (see Eq. \ref{eq_cZ}) estimated from $\fF^{\EMPM}$ with the `ground-truth' image-template correlations measured from $\fF^{\true}$.
For each molecule of interest, three panels are presented in 
Fig. \ref{Fig_viewing_angle_distribution_SET_ONE}. The left-most subplot
is a scatterplot of the image-template correlations. That is, for each 
image $\fA_j$ and each template $(\eazimub_l,\epolara_l)$, we compute a point in two-dimensions with coordinates
\[
(\cZ(\eazimub_l , \epolara_l ; \fA_{j} ; CTF_{j} ; \fF^{\EMPM}),
\cZ(\eazimub_l , \epolara_l ; \fA_{j} ; CTF_{j} ; \fF^{\true})).
\]
With $\nimage = 1024$ images and $N_T = 993$ templates, this yields a data
set with about one million points.
The color in the heatmap corresponds to the log-density of the joint-distribution at that location (see adjacent colorbar).
Note that there is a significant correlation between the estimated- and ground-truth-image-template-correlations.
In the middle subplot for each molecule, 
we show a scatterplot of the template
ranks for each image (see \eqref{templaterankdef}),
rescaled to $[0,1]$ and aggregated over all images.
In the right-most 
subplot, we show a scatterplot of the image ranks for each template
(see \eqref{imagerankdef}),
rescaled to $[0,1]$ and aggregated over all templates.
Once again, the scatterplots are presented as heatmaps, with color corresponding to log-density.

This view of the data sheds light on the differences between maximum-likelihood and maximum-entropy alignment.
Recall that, in standard maximum-likelihood alignment, 
one assigns $\tau_{j}^{\est}$ (for a particular image $\fA_{j}$) to the 
template with the highest template rank.
By contrast, in maximum-entropy alignment, 
one does the reverse: assigning the image with highest image rank to each
sampled template direction.
It is interesting to note that
the estimated and ground-truth image ranks are typically more tightly correlated than the estimated and ground-truth template ranks.
We believe that this phenomenon contributes to the success of our strategy.
Indeed, when the estimated and ground-truth template ranks are highly correlated (as for {\em ISWINCP}), the maximum-likelihood alternating minimization method
(i.e., Bayesian inference with noise-level $0$) succeeds, 
performing roughly as well as our strategy. 
However, when the estimated  and ground-truth template ranks are poorly 
correlated (as for {\em ps1}), then maximum-likelihoood alternating minimization 
does very poorly, often converging to a volume that is worse than random.

\begin{figure}
    \centering
   \includegraphics[width=6.0in]{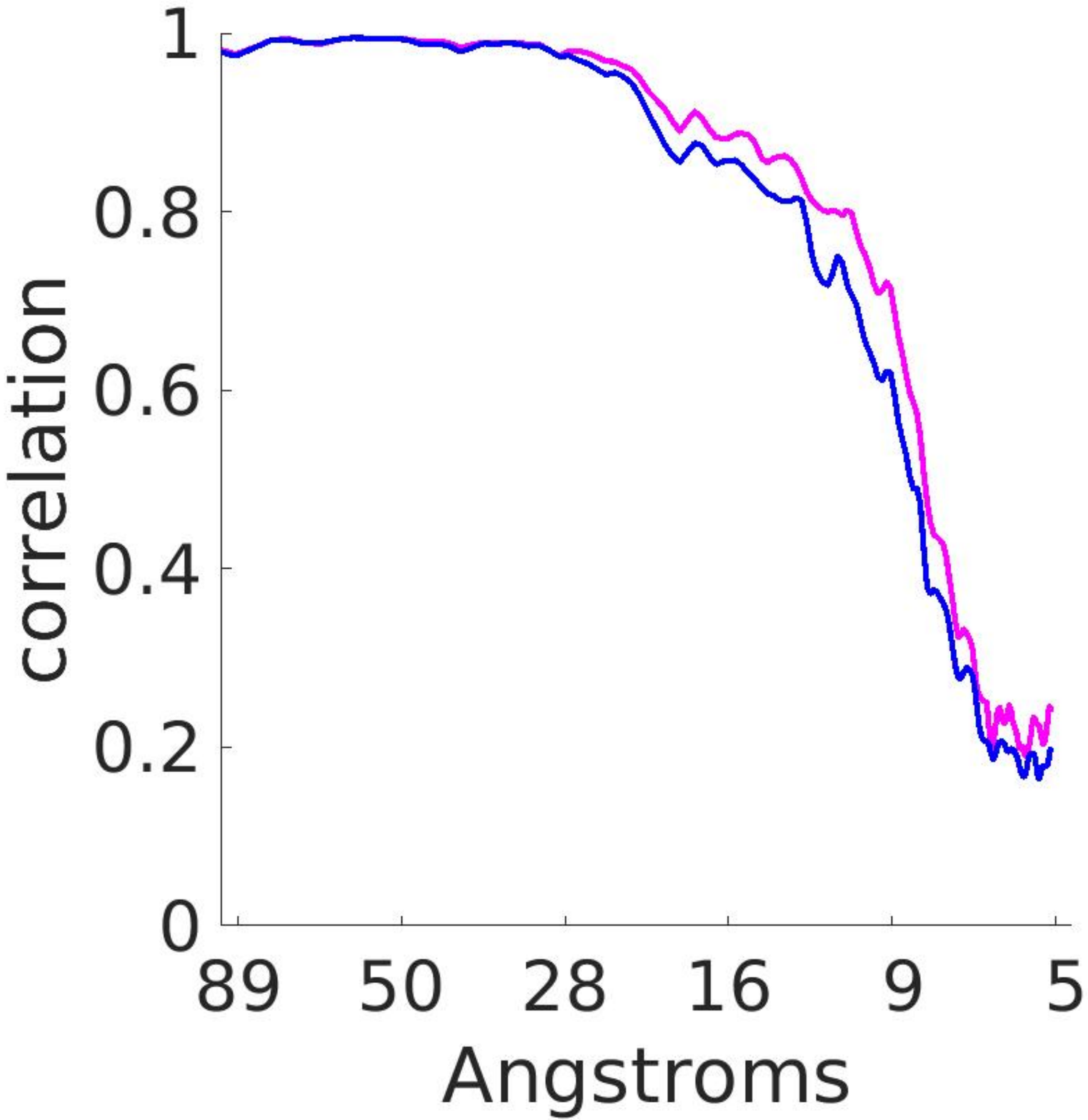} \\
    \caption{
{
    In this figure we show the Fourier shell correlation (FSC) 
associated with a medium-resolution reconstructions produced 
using Relion's {\tt relion\_refine} applied to the first $\nimage=8192$ 
images from the {\em trpv1} dataset (i.e., the first  $\nimage$ picked-particles from the {\tt MRC\_1901} file).
We run this reconstruction twice: once using as initial volume
the RDN reconstruction $\fF^{\RelionDeNovo}$ 
obtained using the first $1024$ images 
and once using as initial volume the EMPM reconstruction $\fF^{\EMPM}$.
We plot the FSC curve obtained by comparing the reconstructed volume 
with the published ground-truth $\fF^{\true}$.
The blue (lower) curve is from the RDN initialization and the 
magenta (higher) curve is from the EMPM initialization.
}
    \label{Fig_trpv1_medres}}
\end{figure}

The illustrations in Figs. \ref{Fig_vol_SET_ONE} and \ref{Fig_viewing_angle_distribution_SET_ONE} are obtained using just one trial of $\fF^{\EMPM}$.
We have carried out extensive simulations and 
these results are fairly typical, with qualitatively similar results with 
repeated runs.
For illustration, we show the values of the volumetric correlations 
collected over multiple trials in Fig. \ref{Fig_collect_SET_ONE}.
The values of $\cZ^{\masked}(\fF^{\EMPM},\fF^{\true})$ shown in magenta are taken over a range of values of $\epsilon_{\tol}\in[0.001,0.010]$, and $\vdlocal\in[0.01,0.10]$.
We did not see a large correlation between either $\epsilon_{\tol}$ or $\vdlocal$ and $\cZ^{\masked}$, provided that $\epsilon_{\tol}$ is sufficiently small and $\vdlocal$ is sufficiently large (i.e., $\epsilon_{\tol}\leq 0.01$ and $\vdlocal\geq 0.01$).
Note that, for these data sets, our strategy reliably produces reconstructions that are similar in quality to $\fF^{\OracleAMStable}$, 
\revTwo{and superior to the results obtained with existing open-source 
packages (DRG, ASP, RDN).}
Finally, we remark that a higher-quality low-resolution reconstruction can 
help improve the results of subsequent refinement.
As an example, we use Relion's {\tt relion\_refine} to produce a medium-resolution reconstruction from the first $8192$ images of this dataset using
two different initial-volumes.
First, we set the initial volume to Relion's de-novo reconstruction 
$\fF^{\RelionDeNovo}$ (referenced in blue in Fig. \ref{Fig_collect_SET_ONE}).
Second, we set the initial volume to the EMPM reconstruction $\fF^{\EMPM}$ shown in Fig. \ref{Fig_vol_SET_ONE}.
Fig. \ref{Fig_trpv1_medres} illustrates the Fourier shell correlations (FSCs) between the resulting medium-resolution reconstructions and the published $\fF^{\true}$.
Note that the higher-quality initial volume $\fF^{\EMPM}$ results in a better medium-resolution reconstruction.
Moreover, as we shall discuss below, better 
low-resolution reconstructions can be useful for other quality-control 
tasks in the cryo-EM pipeline.

We have focused here on a subset of the proteins listed above.
Addition results are deferred to the Supporting material.

{
\begin{remark}
Of note is that the EMPM method works reasonably well even when the particle orientations
are highly non-uniform (see Fig. \ref{Fig_view_angle}).
One might expect that the entropy-maximization
step (which forces images to be assigned to each viewing angle) would cause
difficulties. For de novo low- to moderate-resolution reconstructions, however, it
appears that EMPM benefits from the prevention of spurious clustering
that can develop in pure AM iteration.
At higher resolution, we do not recommend the use of EMPM, as
the entropy-maximization step would result in a loss of resolution and is no longer
necessary once a good initial model has been created.
\end{remark}
}

\begin{remark}
We have implemented our strategy in Matlab, using standard 
(double precision) libraries.
Even for our relatively naive implementation, the total computation time 
for each of our trials is $30-50$ minutes on a desktop workstation.
By comparison, the de-novo reconstruction method within Relion 
typically requires $20-25$ hours on the same machine, while a trial of 
Bayesian inference (also in Matlab) requires $36-48$ hours.
It is also worth noting that,
when constructing $\fF^{\BayesianInference}$, we chose the optimal noise-estimate $\sigma_{\BayesianInference}^{\opt}$. Scanning over multiple values of 
$\sigma_{\BayesianInference}$ would require additional computation time.
In summary, we believe that our strategy has the potential to be extremely 
efficient, especially when properly optimized.
\end{remark}

\subsection{On-the-fly particle rejection}

\begin{figure}
    \centering
    \includegraphics[width=5.0in]{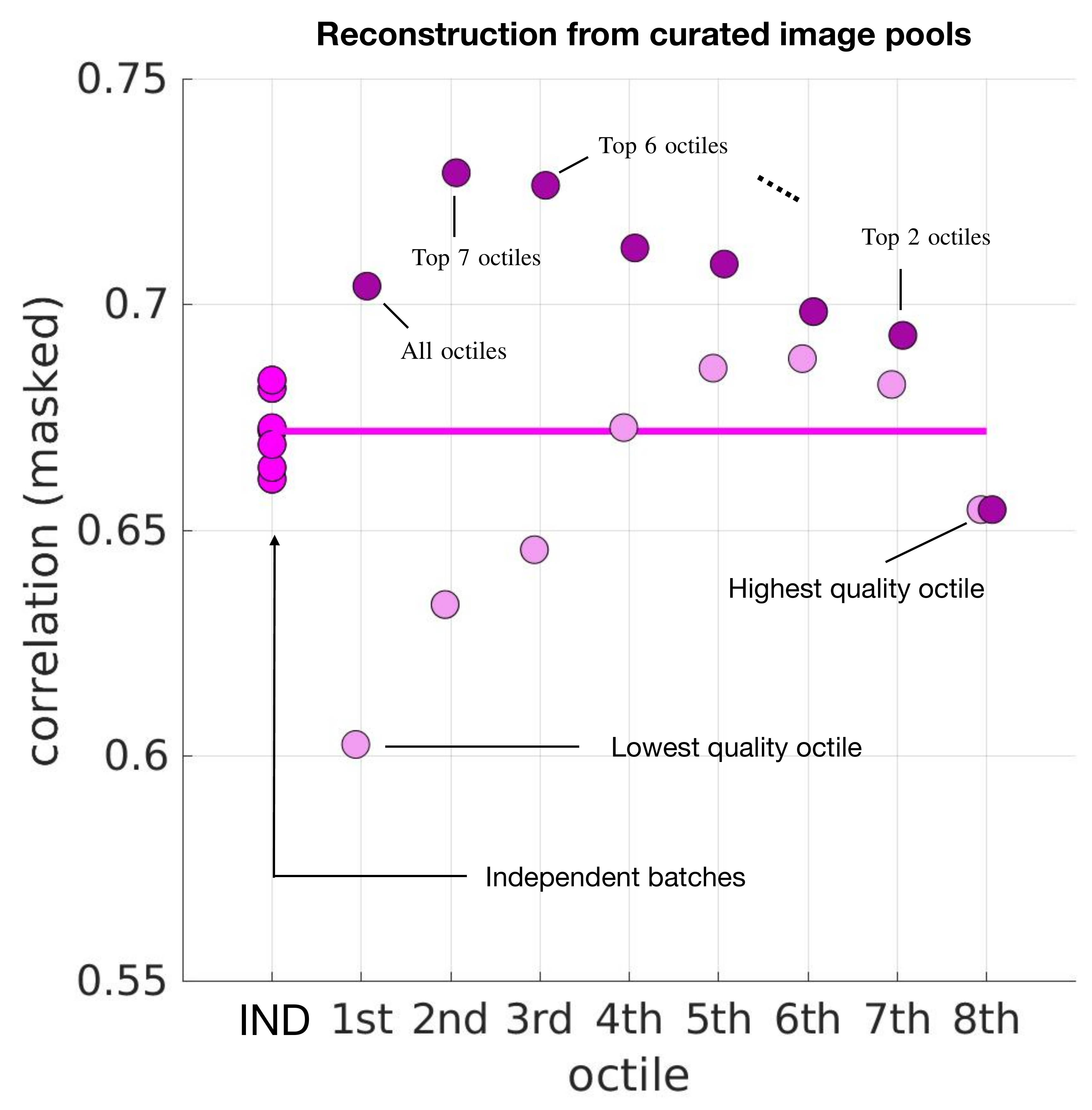}
    \caption{
    In this figure we illustrate the potential for improving reconstruction quality by curating the image-pool.
    For this example we use the {\em ps1} molecule, and consider the first $1024\cdot 8=8192$ images.
    First we divide these images into batches of $\nimage=1024$, and produce multiple reconstructions (magenta, leftmost column).
    The average $\cZ^{\masked}$ for these reconstructions is $\sim 0.67$ (magenta line).
    Next, we use the individual estimates of image-quality (see text) to divide the images into 8 `octiles'.
    We perform a reconstruction using the images within each octile (light-pink).
    Note that those images of poorest quality (i.e., in the first few octiles) produce an inferior reconstruction.
    For this image-pool the images with the very best quality (i.e., in the final octile) tend to be attributed to a small subset of viewing-angles; a reconstruction using only those images is also of rather poor quality.
    We then perform a reconstruction using all the images {\em except} those below a particular octile (dark purple).
    Thus, the reconstruction using all $8192$ images is shown in dark-purple within the `1st' column.
    The reconstruction using the $8192-1024=7168$ images which exclude those in the first octile is shown in dark-purple within the `2nd' column.
    Note that excluding the worst images improves the reconstruction-quality.
    As the number of excluded images grows the number of images used for the reconstruction shrinks; the dark-purple circles are not much higher than the light-pink circles when the octile-index is large.
    }
    \label{Fig_ps1_x0to7_combine}
\end{figure}

Focusing on the {\em ps1} data set, we now consider whether our algorithm can be used to curate particles from small samples.
As shown in Figs. \ref{Fig_viewing_angle_distribution_SET_ONE} and \ref{Fig_viewing_angle_distribution_SET_TWO}, the estimated image-template correlations $\cZ(\eazimub,\epolara;\fA_{j};CTF_{j};\fF^{\EMPM})$ are often quite strongly correlated with the true image-template correlations $\cZ(\eazimub,\epolara;\fA_{j};CTF_{j};\fF^{\true})$.
Motivated by this observation, we define:
\[ \cQ_{j}^{\EMPM} = \max_{t} \cZ(\eazimub_{t},\epolara_{t};\fA_{j};CTF_{j};\fF^\EMPM), \quad 
\cQ_{j}^{\true} = \max_{t} \cZ(\eazimub_{t},\epolara_{t};\fA_{j};CTF_{j};\fF^\true) \]
to be the maximum image-template correlation for each image.
We can use $\cQ_{j}^{\EMPM}$ as a proxy for the `true image-quality' $\cQ_{j}^{\true}$.
We see from Figs. \ref{Fig_viewing_angle_distribution_SET_ONE} and \ref{Fig_viewing_angle_distribution_SET_TWO} that $\cQ_{j}$ typically lies in the interval $[0.1,0.4]$.
It is reasonable to conjecture that
those images with a low estimated quality do not align well to the reconstructed molecule and could be discarded 
when attempting to reconstruct the molecule.
To investigate this hypothesis, we perform a numerical experiment 
with the first $1024\cdot 8=8192$ images in the EMPIAR data set, 
dividing these images into $8$ batches of $\nimage=1024$.
For each batch we produce a reconstruction, using our strategy above.
The $\cZ^{\masked}$ for these single-batch runs are shown as magenta-dots in the `IND' column of Fig. \ref{Fig_ps1_x0to7_combine}.
The average $\cZ^{\masked}$ for these single-batch runs is $\sim 0.67$ (magenta line), with relatively little variation.

Next, we use the estimated image-quality $\cQ_{j}$ to sort the images into octiles (each with $1024$ images). 
Those with the lowest quality scores define the first octile and 
those with the highest quality scores define the eighth octile.
We now perform a reconstruction for each octile separately, with the $\cZ^{\masked}$ shown in light-pink-dots.
Note that the first few octiles produce reconstructions that are quite a bit worse than the typical reconstructions shown in the 1K column.
It is interesting to note that, for this image pool, 
using exclusively the highest quality octile results 
in a poor reconstruction,
since the highest quality images are concentrated around a small subset of viewing-angles.
The result obtained using 
all $\nimage=8192$ images together is shown as the dark-purple dot in 
the `1st' column.
We then show successive
reconstructions using different quality cut-offs. 
Using the top seven octiles ($7168$ images), we 
obtain the improved $\cZ^{\masked}$ shown as a dark-purple dot 
in the `2nd' column.  Using the top six octiles, we get the 
dark purple dot in the `3rd' column, and so forth.
In short, excluding the worst images improves the reconstruction quality, but excluding {\em too many} images causes the quality to drop, 
presumably as a result of losing the coverage from noisier projection directions.

\subsection{Molecules with poor performance on small data sets}

While our results can be considered moderately successful for many of the molecules, our reconstructions for {\em LetB1}, {\em LSUbl17dep} and {\em MlaFEDB} using $1024$
particles are typically quite mediocre.
For {\em LetB1} and {\em LSUbl17dep}, different runs of EMPM yielded 
reconstructions that were quite variable.
We suspect that this inconsistency stems from some underlying 
heterogeneity in the data set itself (see \cite{DTCPLW16}).
For {\em MlaFEDB}, on the other hand, different
trials of EMPM yielded reconstructions $\fF^{\EMPM}$ that were not
close to $\fF^{\true}$, but were often quite close to one another.
This consistency across trials suggests that the structures observed 
in the typical reconstruction might actually exist in the image pool 
we are using. A recurring motif in the reconstructions is that of a 
cylindrical shape with a disjoint structure off to one side.
In this case, we suspect that the image pool contains artifacts 
that are at least 
partially responsible for the mediocre performance of EMPM 
(as well as BI, RDN, DRG, ASP) for this molecule.

\section{Discussion}

As demonstrated above, the EMPM iteration (a modification of alternating-minimization) can produce a reasonable low-resolution reconstruction from a 
small set of raw images.
\revTwo{As seen in Figs. \ref{Fig_collect_SET_ONE} and \ref{Fig_collect_SET_TWO},
it is more robust than many of the currently used 
de-novo molecular reconstruction techniques}, 
and at least as accurate as
Bayesian inference with an optimally chosen estimated noise level.
We did not make an exhaustive comparison with all available packages,
since our main purpose here is to provide a new approach which
is {\em only relevant in the low resolution regime}. 
We believe that the maximum entropy step of EMPM accounts 
for its success, effectively preventing spurious clustering 
of assigned viewing-angles. It should also be noted that the principal mode reduction
described here can be used to accelerate image alignment and volume registration 
in more sophisticated and higher resolution reconstruction algorithms.

EMPM also produces useful estimates of image quality and viewing angle 
distribution -- even when the image pool only contains a few thousand 
images.  These estimates could play a role in online quality-control: 
identifying poor images to be discarded or possibly terminating imaging 
sessions that seem to be yielding poor quality data.

{
Because EMPM is efficient 
(our MATLAB implementation requires less than 5\% of the computation time of RDN and Bayesian inference), it should make cross-validation more {compelling}. 
Our reconstruction requires about 40 min., RDN requires about 20 hrs, full
Bayesian inference requires 40 hrs, cryoDRGN requires 5 hrs and Aspire (the fastest) requires
15 min. (all on the same laptop).
Strong correlation between volumes reconstructed from different subsets of images can 
reasonably be interpreted as  an indication that the data set is of high quality.
On a more speculative note, the image curation presented here uses only image-template
correlations. It is possible that image cross-correlations and the cross-correlations
of the residuals of the least-squares reconstruction step could
also be valuable in settings involving discrete or continuous heterogeneity.}

{
Because EMPM does not require any fine-tuning of parameters, it should be
straightforward to incorporate into existing, widely-used, high-quality 
software pipelines and we hope to do so
in the near future. The code is available from \cite{rangancode}.
}

\newpage

\appendix

\section{Supporting Information} \label{sec:notation}

Various aspects of the EMPM iteration are described
here in greater detail. We begin with 
the discretization of images and volumes
in the Fourier domain, and the techniques used for image alignment
and volumetric reconstruction.
Many of the techniques used require some specification of a `global
tolerance' $\epsilon_{\tol}$, which could be considered a
user-specified parameter. In practice, it appears to 
be sufficient to set it to $10^{-2}$, corresponding
to two digits of accuracy.
In our experiments, we actually scan over $\epsilon_{\tol}\in[0.001,0.01]$ but find that
the quality of the reconstruction does not improve markedly for smaller values of $\epsilon_{\tol}$. In short, EMPM
is not limited by the global tolerance, but presumably
by the experimental error and the structural features of the nonconvex landscape over which we are trying to optimize.

\subsection{Manipulating images in the Fourier domain}
\label{sec_Image_notation}

We use $\vx,\vk \in\Real^{2}$ to represent spatial position and frequency, respectively.
In polar-coordinates these vectors can be represented as:
\begin{eqnarray}
\vx &=& (x \cos \theta, x \sin \theta) \\
\vk &=& (k \cos \psi, k \sin \psi) \text{,}
\end{eqnarray}
for appropriately chosen angles $\theta$ and $\psi$.

\begin{definition}

The Fourier transform of a two-dimensional function 
$A \in \Leb^{2}(\Real^{2})$ is defined as
\begin{equation}
  \fA(\vk) := \iint_{\Real^{2}} A(\vx) \euler^{-\imunit \vk \cdot \vx} \mathop{d\vx}~.
  \label{eq_Ahat}
\end{equation}
\end{definition}

We recover $A$ from $\fA$ using the inverse Fourier transform:
\begin{equation}
A(\vx) = \frac{1}{(2\pi)^{2}} \iint_{\Real^{2}} \fA(\vk) \euler^{+\imunit \vk \cdot \vx} \mathop{d\vk} \text{.}
\end{equation}
The inner-product between two functions $A, B \in \Leb^{2}(\Real^{2})$ is written as
\begin{equation}
  \langle A, B\rangle = \iint_{\Real^{2}} A(\vx)^{\dagger} B(\vx) \mathop{d\vx}
  ~,
\end{equation}
where $z^{\dagger}$ is the complex conjugate of $z \in \Complex$.
We will also use Plancherel's theorem \cite{bracewell},
\begin{equation}
  \langle A, B\rangle = \frac{1}{(2\pi)^{2}} \langle\fA,\fB\rangle
  ~, \qquad \forall A, B \in \Leb^{2}(\Real^{2})~\text{.}
  \label{eq_plancherel}
\end{equation}
We will ignore the factors of $2\pi$ associated with the Fourier transform when they are not relevant.

We represent any given picked-particle image as a function $A \in \Leb^{2}(\Real^{2})$, with values corresponding to the image intensity at each location.
As a consequence of Plancherel's theorem, any inner-product between $A$ and another image (or 2D function) $B$ can be calculated in either real or 
frequency space, the latter of which is often more convenient when aligning images to projections of a volume \cite{Zhao2014,Barnett2017,RSAB20}.

Abusing notation slightly, we will refer to $A(\vk)$ and $\fA(\vk)$ in polar coordinates as:
\begin{eqnarray}
A(x, \theta) &:=& A(\vx) = A(x \cos \theta, x \sin \theta) \\
\fA(k, \psi) &:=& \fA(\vk) = \fA(k \cos \psi, k \sin \psi) \text{.}
\end{eqnarray}
With this notation each $\fA(k,\psi)$ for fixed $k$ and $\psi\in[0,2\pi)$ corresponds to a `ring' in frequency space with radius $k$.

\subsection{Rotation and translation}

Using the notation above, a rotation $\rotation(\egammaz)$ of an image $A$ by angle $\egammaz\in[0,2\pi)$ can be represented as:
\begin{equation}
\label{eq_rotation_definition}
\rotation(\egammaz) \circ A(x, \theta) := A(x, \theta-\egammaz) \text{.}
\end{equation}
Since rotation commutes with the Fourier transform, we have:
\begin{equation}
\wfourier{\left[\rotation(\egammaz) \circ A\right]}(k, \psi) = \rotation(\egammaz) \circ \fA (k, \psi) = \fA(k, \psi-\egammaz) \text{.}
\end{equation}
In this manner, a rotation of any image by $+\egammaz$ can be represented as an angular-shift of each image-ring by $\psi\rightarrow\psi-\egammaz$.

Likewise, a translation $\translation(\vd)$ of an image $A$ by the shift vector $\vd \in \Real^2$ can be represented as:
\begin{equation}
\translation(\vd)\circ A(\vx) := A(\vx - \vd) \text{.}
\end{equation}
In the Fourier domain, this action 
can be expressed as 
\begin{equation}
 \htranslation(\vd,\vk) \odot \fA(\vk) \text{,}
  \label{Tdeltahat}
\end{equation}
where the translation kernel $\htranslation(\vd,\vk)$ is 
\begin{equation}
  \htranslation(\vd,\vk) := \euler^{-\imunit \vd \cdot \vk} \text{.}
  \label{Fdk}
\end{equation}

\subsection{Fourier-Bessel coefficients}

Over the course of our molecular reconstructions, we calculate many 
image-image inner products (see section \ref{sec:correlations}).
To ease the computational burden associated with these calculations, 
it is convenient to use the Fourier-Bessel basis 
\cite{Zhao2014,Zhao2016,Barnett2017,RSAB20}.

To define the Fourier-Bessel coefficients of an image, recall that each image-ring $\fA(k, \psi)$ (for fixed $k$) is a $2\pi$-periodic function of $\psi$, which can itself be represented as a Fourier series in $\psi$:
\begin{equation}
\fA(k, \psi) = \sum_{q=-\infty}^{+\infty} \bA(k;q) \euler^{\imunit q \psi} \text{,}
\end{equation}
for $q\in\Integer$.
The Fourier-Bessel coefficients $\bA(k;q)$ of the image ring $\fA(k,\psi)$ are given by
\begin{equation}
\bA(k;q) = \frac{1}{2\pi} \int_{0}^{2\pi} \fA(k, \psi) \euler^{-\imunit q \psi} \mathop{d\psi}~.
\label{eq_aqk}
\end{equation}
These coefficients can be represented in a more traditional fashion by recalling that the Bessel function $J_{q}(kx)$ can be written as:
\begin{equation}
J_{q}(kx) = \frac{1}{2\pi}\int_{0}^{2\pi} \euler^{-\imunit kx\cos(\psi+\pi/2)-\imunit q \psi} \mathop{d\psi}~ ,
\end{equation}
which, when combined with the definition of the Fourier transform, 
implies that:
\begin{eqnarray}
\bA(k;q) & = & \iint A(x, \theta) \frac{1}{2\pi} \int_{0}^{2\pi} \euler^{-\imunit kx\cos(\psi-\theta) - \imunit q \psi} \mathop{d\psi} xdx\mathop{d\theta}~ \\
 & = & \iint A(x, \theta) \euler^{-\imunit q(\theta+\pi/2)} \frac{1}{2\pi} \int_{0}^{2\pi} \euler^{-\imunit kx\cos(\psi+\pi/2) - \imunit q \psi} \mathop{d\psi} xdx\mathop{d\theta}~\\
 & = & \iint A(x, \theta) \euler^{-\imunit q(\theta+\pi/2)}J_{q}(kx) xdx\mathop{d\theta}~,
\end{eqnarray}
which is the inner product between the original image (in real space) and a 
Bessel function.

Note that, using the Fourier-Bessel representation of an image, the rotation $\rotation_{\egammaz}$ can now be represented as:
\begin{eqnarray}
\rotation_{\egammaz}\fA(k, \psi) & = & \fA(k,\psi-\egammaz) \\
                      & = & \sum_{q=-\infty}^{+\infty} \bA(k;q) \euler^{\imunit q (\psi-\egammaz)} \\
                      & = & \sum_{q=-\infty}^{+\infty} \bA(k;q)\cdot\euler^{-\imunit q\egammaz} \euler^{\imunit q \psi} \text{,}
\label{eq_fourier_bessel_rotation}
\end{eqnarray}
such that the Fourier-Bessel coefficients of the rotated image ring $\rotation_{\egammaz}\circ \fA(k,\cdot)$ are given by the original Fourier-Bessel coefficients $\bA(k;q)$, each multiplied by the phase-factor $\euler^{-\imunit q\egammaz}$.
Note that \eqref{eq_fourier_bessel_rotation} easily allows for arbitrary rotations by any $\egammaz$, requiring only a pointwise multiplication by the appropriate phase factor.
This observation allows for the efficient alignment of one image to another (see \cite{Barnett2017,RSAB20}).

\subsection{Image discretization}
\label{sec_Image_discretization}

We denote by $\Omega(1)$ and $\Omega(\kmax)$ the ball of radius $1$ and $\kmax$, respectively, in either real  or frequency space:
\begin{equation}
\Omega(1) := \{ \vx \in \Real^{2} \mbox{~such~that~} \|\vx\| \leq 1 \} \text{,} \quad \Omega(\kmax) := \{ \vk \in \Real^{2} \mbox{~such~that~} \|\vk\| \leq \kmax \} \text{,}
\end{equation}
and we will assume that all the images considered are supported 
in $\Omega(1)$. The support constraint
$A(\vx) \subseteq \Omega(1)$ implies that 
the representation $\fA(\vk)$ will have a bandlimit of $1$, so 
that $\fA$ can be accurately reconstructed from its values sampled on a frequency grid with spacing $\bigO(1)$ \cite{finufft}.

We also assume that any relevant signal within the images has a maximum effective spatial-frequency magnitude of $\kmax$; i.e., that the salient features of $\fA$ are concentrated in $\vk\in\Omega_{\kmax}$.
Consequently, we expect that the inversion
\begin{equation}
A(\vx) \approx \frac{1}{(2\pi)^{2}} \iint_{\Omega_\kmax} \fA(\vk)\euler^{\imunit \vk \cdot \vx} \mathop{d\vk}
\end{equation}
accurately reconstructs $A$ to the resolution associated with $\kmax$
(ignoring high frequency ringing artifacts).
When these assumptions hold we won't lose much accuracy when applying these transformations in a discrete setting (as discussed below).

Note that the largest $\kmax$ that we can reasonably expect to attain is related to the number of pixels `$\npixel$' in the original image.
More specifically, once we assume that $A(\vx)$ is supported in $\vx\in\Omega(1)$, the pixel-spacing will be $\dx = 2/\npixel$, implying the Nyquist spatial-frequency is given by $\pi/\dx = \npixel\pi/2$.

Since we are using $\nimage$ picked-particle images $A_{j}$ to produce a low-resolution estimate of the imaged molecule (denoted $F$ below),
we will typically choose $\kmax$ to be significantly smaller than $\npixel\pi/2$, often in the range of $\kmax\sim 50$.
Given the pixel spacing, this corresponds to a low-resolution reconstruction with approximately $20$\AA resolution.
Additionally, the fact that $x\leq 1$ implies that the Bessel coefficients $\bA(k,q)$ will be concentrated in the range $|q|\lesssim k$, meaning that the Bessel coefficients $\bA(k,q)$ across all $k\in[0,\kmax]$ will be concentrated in $q\in[-\qmax/2,+\qmax/2-1]$ for $\qmax=\bigO(\kmax)=\bigO(\npixel)$.

With the notation above, we can consider an $\npixel\times\npixel$ image as a discrete set of pixel-averaged samples within $[-1,+1]^{2}$:
\begin{equation}
A_{n_{1},n_{2}} = \frac{1}{\dx^{2}}\int_{n_{1}\dx}^{(n_{1}+1)\dx} \int_{n_{2}\dx}^{(n_{2}+1)\dx} A(x_1,x_2) \, dx_2 dx_1
\end{equation}
for indices $n_{1},n_{2}\in \{0,\ldots,\npixel-1\}$.
We approximate the Fourier transform $\fA$ at any $\vk \in \Real^{2}$ via the simple summation:
\begin{equation}
\fA(\vk) = (\dx)^{2} \sum_{n_{2}=0}^{\npixel-1} \sum_{n_{1}=0}^{\npixel-1} A_{n_{1},n_{2}} e^{-\imunit \vk \cdot \vx_{n_{1},n_{2}}}~ \text{,}
\label{eq_Ahattrap}
\end{equation}
where $\vx_{n_{1},n_{2}}$ is the appropriately-chosen pixel-center $\dx\left(n_{1}+\frac{1}{2},n_{2}+\frac{1}{2}\right)$.
Because we have assumed that the image is sufficiently well sampled (i.e., that $\fA$ contains little relevant frequency-content above the Nyquist-frequency $\kmax$), we expect the simple sum above to be accurate.

We will typically evaluate $\fourier{A}(\vk)$ for $\vk$ on a polar-grid, with $k$- and $\psi$-values corresponding to suitable quadrature nodes, 
using the non-uniform FFT (NUFFT) to compute $\fA(k,\psi)$ at those nodes
(see \cite{finufft,RSAB20}).
The quadrature in $k$ corresponds to a set of $\rmax$ nodes: 
$k_{1},\ldots,k_{\rmax}$ and radial weights $w_{1},\ldots,w_{\rmax}$ so that
\begin{equation}
\int_{0}^{\kmax} g(k)kdk \approx \sum_{r=1}^{\rmax} g(k_{r})w_{r} 
\end{equation}
to high accuracy for a smooth function $g(k)$.
The exact radial quadrature nodes we select are usually those inherited from a volumetric grid (see section \ref{sec_Volume_notation} below).
For convenience, we will sometimes refer to the weight function
$w(k)$ as some continuous interpolant of the quadrature-weights, 
with $w(k_{r})=w_{r}$.
The $\qmax$ angular-nodes $\psi_{0},\ldots,\psi_{\qmax-1}$ are equispaced in the periodic interval $[0,2\pi)$, with a spacing of $\dpsi=2\pi/\qmax$, and $\psi_{q'}=q'\dpsi$.
These equispaced $\psi$-nodes allow for spectrally accurate trapezoidal quadrature in the $\psi$-direction, and we approximate the Fourier-Bessel coefficients of each image-ring $\fA(k,\psi)$ as follows:
\begin{equation}
\bA(k,q) \approx \sum_{q^{\prime}=0}^{\qmax-1} \fA(k,\psi_{q^{\prime}}) \exp\left(-\imunit q \psi_{q^{\prime}}\right) \dpsi \text{,}
\end{equation}
with the index $q$ considered periodically in the interval $[-\qmax/2,\ldots,+\qmax/2-1]$ (so that, e.g., the $q$-value of $\qmax-1$ corresponds to the $q$-value of $-1$).

\subsection{Inner products}

The inner-product between an image $A$ and a rotated and translated version of image $B$ is denoted by:
\begin{eqnarray}
\cX(\egammaz , \vd ; A,B) &:=& \langle \rotation(\egammaz)\circ A , \translation(\vd)\circ B \rangle \text{,}
\end{eqnarray}
and (ignoring factors of $2 \pi$) is approximated by:
\begin{eqnarray}
\cX(\egammaz , \vd ; \fA,\fB) & \approx & \iint_{\Omega(\kmax)} \left[\rotation(\egammaz)\circ\fA(\vk)\right]^{\dagger} \cdot \left[\htranslation(\vd,\vk)\odot\fB(\vk)\right] d\vk  \text{.}
  \label{Xdg}
\end{eqnarray}
The expression $\cX$ can be thought of as a {\em bandlimited} inner product.

\subsection{Volume notation}
\label{sec_Volume_notation}

We use $\vx,\vk \in\Real^{3}$ to represent spatial position and frequency, respectively.
In spherical coordinates, the vector $\vk$ is represented as:
\begin{eqnarray}
\vk &=& k\cdot \hk, \text{\ \ with\ \ } \hk = (\cos \kazimub \sin \kpolara, \sin \kazimub \sin \kpolara , \cos \kpolara) \text{,}
\end{eqnarray}
with polar angle $\kpolara$ and azimuthal angle $\kazimub$ representing the unit vector $\hk$ on the surface of the sphere $S^{2}$.

Using a right-handed basis, a rotation about the third axis by angle $\eazimub$ is represented as:
\begin{equation}
\rotation_{z}(\eazimub) =
\left(
\begin{array}{ccc}
+\cos\eazimub & -\sin\eazimub & 0 \\
+\sin\eazimub & +\cos\eazimub & 0 \\
0            & 0            & 1
\end{array}
\right) \text{,}
\end{equation}
and a rotation about the second axis by angle $\epolara$ is represented as:
\begin{equation}
\rotation_{y}(\epolara) =
\left(
\begin{array}{ccc}
+\cos\epolara & 0 & +\sin\epolara \\
0            & 1 & 0            \\
-\sin\epolara & 0 & +\cos\epolara 
\end{array}
\right) \text{.}
\end{equation}

A rotation $\rotation(\tau)$ of a vector $\vk\in\Real^{3}$ can be represented by the vector of Euler angles $\tau = (\eazimub, \epolara, \egammaz)$:
\begin{equation}
\rotation(\tau)\cdot\vk = \rotation_{z}(\eazimub)\circ \rotation_{y}(\epolara)\circ \rotation_{z}(\egammaz)\cdot \vk \text{.}
\end{equation}

We represent a given volume as a function 
$F\in\Leb^{2}(\Real^{3})$, with values corresponding to the intensity 
at each location $\vx\in\Omega(1)$.
We will refer to $\fF(\vk)$ in spherical coordinates as $\fF(k, \hk)$.
With this notation, each $\fF(k,\cdot)$ corresponds to a `shell' in frequency-space with radius $k$.
The rotation of any volume $\rotation(\tau)\circ\fF(k,\hk)$ corresponds to the function $\fF(k,\rotation(\tau^{-1})\cdot\hk)$, where $\tau^{-1}$ corresponds to the vector of Euler angles 
$\tau^{-1} = \left(-\egammaz,-\epolara,-\eazimub \right)$.

\subsection{Spherical harmonics}

Using the notation above, we can represent a volume $\fF(k,\hk)$ as:
\begin{equation}
\fF(k,\hk) = \sum_{l=0}^{+\infty}\sum_{m=-l}^{m=+l} \sF(k;l,m) Y_{l}^{m}(\hk) \text{,}
\end{equation}
where $Y_{l}^{m}(\hk)$ represents the spherical harmonic of degree $l$ 
and order $m$:
\begin{equation}
Y_{l}^{m}(\kpolara,\kazimub) =  Z_{l}^{m} \euler^{+\imunit m \kazimub}P_{l}^{m}(\cos\kpolara) \text{,}
\end{equation}
with $P_{l}^{m}$ representing the usual (unnormalized) associated Legendre 
function, and $Z_{l}^{m}$ the normalization constant:
\begin{equation}
Z_{l}^{m} = \sqrt{\frac{2l+1}{4\pi}\times\frac{\left(l-|m|\right)!}{\left(l+|m|\right)!}} \text{.}
\end{equation}
The coefficients $\sF(k;l,m)$ define the spherical harmonic expansion of the $k$-shell $\fF(k,\cdot)$.

Given a rotation $\tau=(\eazimub,\epolara,\egammaz)$, 
we can represent the rotated volume $\fG:=\rotation(\tau)\fF$ as:
\begin{equation}
\sG_{l}^{m_{1}} = \sum_{m_{2}=-l}^{m_{2}=+l} \euler^{-\imunit m_{1} \egammaz}d_{m_{1},m_{2}}^{l}(\epolara) \euler^{-\imunit m_{2} \eazimub} \sF_{l}^{m_{2}} \text{,}
\end{equation}
where $d_{m_{1},m_{2}}^{l}(\epolara)$ represents the 
Wigner d-matrix of degree $l$ associated with the second 
Euler angle $\epolara$.

\subsection{Volume discretization}
\label{sec_Volume_discretization}

As with our discretization of images, we assume that the volume $F(\vx)$ is supported in $\Omega(1)$, and that the relevant frequency content 
is contained in $\Omega(\kmax)$.

We then discretize the radial component of $\Omega(\kmax)\in\Real^{3}$ using a Gauss-Jacobi quadrature for $k$ built with a weight-function corresponding to a radial weighting of $k^{2}dk$; the number of radial quadrature nodes $\rmax$ will be of the order  $\bigO(\kmax)$.
On each of the $k$-shells, we build a quadrature for $\hk$ that discretizes the polar-angle $\epolara$ using a legendre-quadrature on $[0,\pi]$, and the azimuthal-angle $\eazimub$ using a uniform periodic grid in $[0,2\pi)$, with associated weights designed for integration on the surface of the sphere; the total number of spherical-shell-quadrature-nodes $\tmax(k_{r})$ for any particular $k_{r}$ need only be $\bigO(k_{r}^{2})$ \cite{Barnett2017}.

Each of the $k$-shells $\fF(k_{r},\cdot)$ can be accurately described using spherical harmonics of order $l=\bigO(k_{r})$.
Thus, the number of spherical-harmonic coefficients required for each shell $\fF(k_{r},\cdot)$ is $\bigO(k_{r}^{2})$.
The total number of spherical-harmonic coefficients required to approximate $\fF(k,\hk)$ over $\Omega(\kmax)$ is $\bigO(\kmax^{3})$, with a maximum degree of $\lmax=\bigO(\kmax)$.
The associated maximum order will be $\mmax=1+2\lmax$, which is also $\bigO(\kmax)$.

\subsection{Template generation}

Given any volume $\fF(\vk)$, we can generate templates efficiently using the techniques described in \cite{Barnett2017}, which make use of the spherical harmonic representation $\sF(k;l,m)$.
For our purposes, we generate the templates described in Eq. \ref{eq_translated_CTF_modulated_template} for a collection of viewing directions
$(\eazimub_{t},\epolara_{t})$, for $t\in\{1,\ldots,\tmax\}$.
The particular choice of the third Euler angle ($\egammaz_{t}$) for each of 
the templates is dictated by the strategy for template generation 
in $\cite{Barnett2017}$, and does not play an important role,
since it corresponds to an in-plane rotation of the relevant slice.
We discretize the set of viewing directions
to correspond with the $\hk$ discretization on the largest $k$-shell in 
our volumetric (spherical) quadrature grid:
that is, $\tmax:=\tmax(k_{\rmax})$.
These viewing directions will be used to enumerate the templates (i.e., $\ntemplate=\tmax$).
Each of the templates are stored on a polar or Fourier-Bessel grid 
identical to those used to store the images.

\subsection{Image noise model}
\label{sec_Image_noise_model}

We assume that each of the picked-particle images $\fA_{j}$ is a noisy version of a `signal' corresponding to a 2-dimensional projection of some (unknown) molecular electron-density-function $F^{\true}$:
\begin{equation}
\fA^{\signal}(\vk) = 
    \translation(+\vd_{j}^{\true},\vk)\odot
\fS(\tau_{j}^{\true} ; CTF ; \fF)(\vk)\text{,}
\end{equation}
where $\tau_{j}^{\true}$ and $\vd_{j}^{\true}$ correspond to the true (but unknown) viewing-angle and displacement associated with image $\fA_{j}$.
In this context we say the image's `true viewing direction' is defined by the first two Euler angles in $\tau_{j}^{\true}$, namely the azimuthal and 
polar angles $\eazimub_{j}^{\true}$ and $\epolara_{j}^{\true}$.

We will model the noise in real space as independent and identically 
distributed ({\em iid}), with a variance of $\sigma^{2}$ on a unit scale in 
two-dimensional real space.
This is a simple model of detector noise 
\cite{Sigworth1998,Scheres2009,Sigworth2010,Lyumkis2013}, and does not 
take into account more complicated features, such as structural noise 
associated with image preparation (see \cite{Baxter2009}).
We assume that the image values are the sum 
of the signal plus noise:
\begin{equation}
A_{n_{1},n_{2}} = A(\vx_{n_{1},n_{2}}) = A^{\signal}_{n_{1},n_{2}} + A^{\noise}_{n_{1},n_{2}} \text{,} 
\end{equation}
where the signal and noise are represented by the arrays $A^{\signal}_{n_{1},n_{2}}$ and $A^{\noise}_{n_{1},n_{2}}$, respectively.
Because each entry in $A^{\noise}_{n_{1},n_{2}}$ corresponds to an average over the area-element $\dx^{2}$ associated with a single pixel, we expect the variance of each $A^{\noise}_{n_{1},n_{2}}$ to scale inversely with $\dx^{2}$.
That is to say, we'll assume that each element of the noise array $A^{\noise}_{n_{1},n_{2}}$ is drawn from the standard normal distribution $\cN\left(0,\frac{\sigma^{2}}{\dx^{2}}\right)$, with zero mean and a variance of $\sigma^{2}/\dx^{2}$.
With these assumptions, downsampling the image by averaging neighboring pixels will correspond to an increase in the effective pixel-size and a simultaneous reduction in the variance of the noise associated with each (now larger) pixel.

Given the assumptions above, $\fA$ will be modeled by:
\begin{equation}
\fA(\vk) = \fA^{\signal}(\vk) + \fA^{\noise}(\vk) \text{,} 
\end{equation}
where $\fA^{\noise}$ is now complex and {\em iid}.
Because the noise-term $A^{\noise}$ is real, the noise-term $\fA^{\noise}$ will be complex, with the conjugacy constraint that $\fA^{\noise}(+\vk) = \fA^{\noise}(-\vk)^{\dagger}$.

\subsection{Variance Scaling}
\label{sec_Variance_Scaling}

We define $\fsigma^{2}=\frac{\pi^{2}\sigma^{2}}{\dx^{2}}$ as the 
variance on a unit scale in two-dimensional frequency space.
We expect that the noise term $\fA^{\noise}(\vk)$ integrated over any 
area element $\dk^{2}$ will have a variance of 
$\frac{\fsigma^{2}}{\dk^{2}}$.
Recalling our polar quadrature described above, we typically sample 
$\vk$ along $\rmax$ radial quadrature nodes $k_{1},\ldots,k_{\rmax}$ 
and $\qmax$ angular quadrature nodes $\psi_{0},\ldots,\psi_{\qmax-1}$, 
with radial and angular weights $w_{r}$ and $\dpsi$, respectively.
In this quadrature scheme, each quadrature node 
$\vk_{rq}=(k_{r},\psi_{q})$ is associated with the area element 
$w_{r}\dpsi$, which approximates the `$kdk$' integration weight.
Consequently, we expect that the noise term $\fA^{\noise}(\vk_{rq})$ 
evaluated at $(k_{r},\psi_{q})$ on our polar quadrature grid will have 
a variance of $\frac{\fsigma^{2}}{w_{r}\dpsi}$.

Given the radially-dependent variance associated with each degree of freedom in the images, we will typically consider rescaled images $\fA(k_{r},\psi_{q})\sqrt{w_{r}\dpsi}$, rather than the raw image values $\fA(k_{r},\psi_{q})$.
Fro simplicity, we will often write this as $\fA(\vk)\sqrt{w(k)}$, where we assume that $\dpsi$ is the same for every $k$-value, and $w(k)$ is some 
functional interpolant of the quadrature weights with $w(k_{r})=w_{r}$.
This rescaling allows us to equate the standard vector 2-norm with the 
$L^2$ integral (up to quadrature error):
\[ 
\sum_{r=1}^{r=\rmax}\sum_{q=1}^{q=\qmax} \| \fA(k_{r},\psi_{q}) \cdot \sqrt{w_{r}\dpsi} \|^{2} 
\approx
\int_{k=0}^{k=\kmax}\int_{\psi=0}^{\psi=2\pi} \|\fA(k,\psi)\|^{2} kdkd\psi
\text{.}
\]
This rescaling will allow our definition of the radial objective function 
below to correspond to a standard log-likelihood 
(see \eqref{eq_Ccost} in section \ref{sec_Choice_of_principal_modes}).
In a similar fashion, the least-squares formulation we will use for 
volume reconstruction corresponds to a maximum-likelihood estimate, 
with similar properties holding for the inner products used for alignment 
(see sections \ref{sec_Volume_Reconstruction} and \ref{sec_Image_Alignment} below).

\subsection{Radial Principal Modes}
\label{sec_Radial_Principal_Modes}

In these sections, we briefly review the notion of {\em radial principal modes} for image compression, described in \cite{Rangan22}.
The premise is simple: given an image $\fA(k,\psi)$, 
not all the frequency rings $\fA(k,\cdot)$ are equally useful.
By selecting the appropriate linear combinations of frequency rings, 
we can compress the images while retaining information useful for alignment.
We choose the coefficients of this linear combination by solving 
an eigenvalue problem, which is why we refer to these linear combinations as principal modes.

This strategy has two benefits.
The first is to make the overall computation more efficient; rather than 
retaining all $\rmax\approx 50$ radial quadrature nodes, 
we can often compress the majority of relevant information into a smaller 
number of principal modes (typically 10-16).
The second benefit is that, by retaining only those principal modes which 
are useful for alignment, we explicitly avoid those degrees of freedom 
within the original images that are dominated by noise and not useful for 
alignment. The same set of modes can be used for 
the experimental images and the estimated volumes.
We make use of this compression for both the alignment and 
reconstruction steps of our EMPM iteration.

\subsection{Principal image rings and principal volume shells}
\label{sec_Principal_image_rings}

Let us assume that we are given a unit vector
$\vu=\transpose{[u_{1},u_{2},\ldots,u_{\rmax}]}\in\Real^{\rmax}$. 
We define the linear combination $[\transpose{\vu}\fA](\psi_{q})$ as:
\begin{eqnarray}
[\transpose{\vu}\fA](\psi_{q}) = \sum_{r=1}^{\rmax} u_{r}\fA(k_{r},\psi_{q})\sqrt{w_{r}\dpsi} \text{,}
\end{eqnarray}
where the rescaling factor $\sqrt{w_{r}\dpsi}$ is designed to equalize the variances associated with the different $k$-values (see section \ref{sec_Variance_Scaling}).

In a moment, we will choose our $\vu$ to be one of the eigenvectors of 
an $\rmax\times\rmax$ matrix, and we will refer to the linear combination $[\transpose{\vu}\fA](\psi_{q})$ as the `principal image ring' associated with 
$\vu$.
Note that, based on our noise model, we assume that the image noise 
is {\em iid} across image-rings.
Consequently, the variance of the noise $[\transpose{\vu}\fA^{\noise}](\psi_{q})$ is equal to $\sum_{r}\vu_{r}^{2}\fsigma^{2}$, which is equal to $\|\vu\|^{2}\fsigma^{2}$.
Because $\vu$ is a unit vector, this last expression is simply $\fsigma^{2}$, which is the same as the variance of any individual term $\fA^{\noise}(k_{r},\psi_{q})$ for any particular $k_{r}$.
Moreover, if we consider two orthonormal vectors $\vu_{1}$ and $\vu_{2}$, then $[\transpose{\vu_{1}}\fA^{\noise}](\psi_{q})$ and $[\transpose{\vu_{2}}\fA^{\noise}](\psi_{q})$ will be independent random variables, each drawn from $\cN(0,\fsigma^{2})$.

The same strategy can be used to define a principal volume shell 
$[\transpose{\vu}\fF](\hk)$ as:
\begin{eqnarray}
[\transpose{\vu}\fF](\hk) = \sum_{r=1}^{\rmax} u_{r}\fF(k_{r},\hk)\sqrt{w_{r}\dpsi} \text{,}
\end{eqnarray}
using the same weights as those used for the principal image rings.
In most situations, we will assume that the angular discretization 
$\dpsi$ is uniform across all the image rings, allowing us to 
ignore the $\sqrt{\dpsi}$ factor in the rescaling-factor above.

\subsection{Choice of principal modes}
\label{sec_Choice_of_principal_modes}

To select the actual principal modes used above, 
we construct the following objective function for the vector $\vu$:
\begin{eqnarray}
\Ccost(\vu ; F , \mu_{\tau} , \mu_{\vd} , \mu_{CTF} )
& = & 
\int\cdots\int
\left\Vert
\transpose{\vu} \fS(\tau,\vd;CTF,F) - \transpose{\vu} \fS(\tau',\vd';CTF',F)
\right\Vert ^{2}
\cdot \nonumber \\
& & 
\hspace{-2cm}
d\psi
d\mu_{\tau}\left(\tau\right)
d\mu_{\tau}\left(\tau'\right)
d\mu_{\vd}\left(\vd\right)
d\mu_{\vd}\left(\vd'\right)
d\mu_{CTF}(CTF)
d\mu_{CTF}(CTF')
\text{,}
\label{eq_Ccost}
\end{eqnarray}
where 
$\fS(\tau , \vd ; CTF ; F)$ is the 
usual noiseless image template associated with viewing-angle $\tau$, 
displacement $\vd$ and the contrast-transfer-function $CTF(\vk)$.
The integral in \eqref{eq_Ccost} is taken over the distributions of viewing angle, displacement and CTF function 
(i.e., $\mu_{\tau}$, $\mu_{\vd}$ and $\mu_{CTF}$, respectively).
Given a volume $F$, as well as distributions $\mu_{\tau}$, $\mu_{\vd}$ and $\mu_{CTF}$, the objective function $\Ccost(\vu)$ is (up to an affine transformation) equal to the negative log-likelihood of mistaking one noiseless principal image ring (constructed using principal-vector $\vu$) for a randomly rotated version of another (also constructed using $\vu$).
Those vectors $\vu$ for which $\Ccost(\vu)$ is high correspond to linear combinations of frequencies which are useful for alignment (i.e. for which the typical principal image rings look very different from one another).
Conversely, those vectors $\vu$ for which $\Ccost(\vu)$ is low correspond to linear combinations of frequencies which are not useful for alignment (i.e., for which the typical principal image rings look quite similar to one another).

The objective function $\Ccost(\vu)$ is quadratic in $\vu$, which means that it can be optimized by finding the dominant eigenvector of the Hessian $\Ckern := \partial_{\vu\vu}\Ccost$, defined as the $\rmax\times\rmax$ matrix:
\[ \Ccost(\vu) = \transpose{\vu} \cdot \Ckern \cdot \vu \text{.} \]
Moreover, the corresponding eigenvalue will be equal to the 
objective function $\Ccost$ for that eigenvector.
To actually compute the Hessian $\Ckern$ in practice, 
we use one of two strategies.

\subsection{Empirical principal modes}

At first, when we have no good estimate of either $F$ or the 
distributions $\mu_{\tau}$, $\mu_{\vd}$ or $\mu_{CTF}$, we simply 
use the images $A_{j}$ themselves to form a Monte Carlo estimate of the 
integral in \eqref{eq_Ccost}.
This boils down to the following simple approximation:
\begin{eqnarray}
\Ccost(\vu)^{\empirical}
& = &
\frac{1}{\nimage^{2}}\sum_{j,j'}
\iiint
\left\Vert
\transpose{\vu} \rotation(\egammaz)\fA_{j} - \transpose{\vu} \rotation(\egammaz')\fA_{j'}
\right\Vert ^{2}
d\psi d\egammaz d\egammaz'
\text{,}
\label{eq_Ccost_empirical}
\end{eqnarray}
where we have replaced the noiseless images 
$\fS(\tau , \vd ; CTF ; F)$ and $\fS(\tau , \vd ; CTF ; F)$ 
in the integrand of Eq. \ref{eq_Ccost} with randomly selected images 
$\fA_{j}$ and $\fA_{j'}$ from the original image pool, and used the 
empirical distribution of the images themselves to stand in for the 
unkown $\mu_{\tau}$, $\mu_{\vd}$ and $\mu_{CTF}$ 
(after averaging over the in-plane rotations associated with 
$\egammaz$ and $\egammaz'$).
The associated kernel $\Ckern^{\empirical}$ then reduces to:
\begin{equation}
\Ckern(r,r')^{\empirical}  \propto  \sqrt{w_{r}w_{r'}}\cdot 
\left( \frac{1}{\nimage}\sum_{j} \sum_{q} \bA_{j}(k_{r},q)^{\dagger}\bA_{j}(k_{r'},q) - 
\frac{1}{\nimage^2}\left[\sum_{j}\bA(k_{r},0)\right]^{\dagger}\left[\sum_{j}\bA(k_{r'},0)\right] \right) \text{.}
\end{equation}

\subsection{Estimated principal modes}

Once we have a reasonable estimate of the imaged volume $F^{\estimated}$, we can use it (along with the current estimated distributions $\mu_{\tau}^{\estimated}$ and $\mu_{\vd}^{\estimated}$) to calculate $\Ccost^{\estimated}$ via \eqref{eq_Ccost}.
In practice, we typically assume that $\mu_{\tau}^{\estimated}$ is uniform and $\mu_{\vd}^{\estimated}$ is an isotropic Gaussian centered at the origin 
with a standard deviation estimated using the current $\vd_{j}^{\estimated}$.
With these assumptions the associated kernel $\Ckern^{\estimated}$ can be computed analytically from the spherical harmonic coefficients $\sF(k;l,m)$:
\begin{eqnarray}
\Ckern_{r,r'}^{\estimated}
& \propto & 
+\sqrt{W_{r}W_{r'}}\cdot
\left[ 
\sum_{l=0}^{+\infty} \sum_{m=-l}^{m=+l} \sF(k_{r};l,m)^{\dagger}\sF(k_{r'};l,m) 
\right] 
\cdot \boldsymbol{E}^{+}(k_{r},k_{r'}) 
\\
& & 
- 
\sqrt{W_{r}W_{r'}}\cdot
\left[ 
\sF(k_{r};0,0)^{\dagger}\sF(k_{r'};0,0) 
\right] 
\cdot \boldsymbol{E}^{-}(k_{r}) \cdot \boldsymbol{E}^{-}(k_{r'})
\text{,}
\label{eq_template_cost_with_translations}
\end{eqnarray}
as $\qmax\rightarrow\infty$, 
where the CTF-modulated weight $W_{r}$ takes the form:
\begin{eqnarray}
\sqrt{W_{r}} & = & 
\frac{\sqrt{w_{r}}}{\nimage}\sum_{j} CTF_{j}(k_{r}) \text{,}
\end{eqnarray}
and the terms $\boldsymbol{E}^{+}$ and $\boldsymbol{E}^{-}$ denote
\begin{eqnarray}
  \boldsymbol{E}^{+}(k_{r},k_{r'}) = \exp\left(-\frac{\sigma_{\vd}^{2}}{2}\left[ k_{r}^{2} + k_{r'}^{2} \right]\right) \cdot \exp\left( k_{r}k_{r'}\sigma_{\vd}^{2} \right) \text{,}
\end{eqnarray}
and
\begin{eqnarray}
\boldsymbol{E}^{-}(k_{r}) = \tilde{k}_{r}\sqrt{\frac{\pi}{2}} \exp\left(-\tilde{k}_{r}^{2}\right) \left( {\cal I}_{-1/2}(\tilde{k}_{r}^{2}) - {\cal I}_{+1/2}(\tilde{k}_{r}^{2}) \right) \text{,}
\end{eqnarray}
where $\tilde{k}_{r} = k_{r}\sigma_{\vd}/2$, and ${\cal I}_{q}$ refers to the modified Bessel function of first kind of order $q$.

\subsection{Multiple principal modes}

We select the principal modes $\vu_{1},\ldots,\vu_{\nrank}$ to be the dominant $\nrank$ eigenvectors of the kernel $\Ckern$.
We typically choose $\nrank$ such that the $(\nrank+1)$-th eigenvalue 
of $\Ckern$ is less than the global tolerance times the first (dominant) 
eigenvalue of $\Ckern$.
When we set the global tolerance to be $\epsilon_{\tol}\sim 1e-2$ in our numerical experiments below, this selection criterion typically corresponds to 
selecting $12-16$ principal modes from $\Ckern^{\empirical}$, and slightly fewer ($10-14$) from $\Ckern^{\estimated}$.

We will refer to the matrix $U$ as the collection of $\nrank$ 
principal modes:
\begin{equation}
U = \left[ \begin{array}{cccc} u_{1} & u_{2} & \cdots & u_{\nrank} \end{array} \right] \text{,}
\label{eq_U_nrank}
\end{equation}
and refer to the collection of respective principal images as $\transpose{U}\fA$:
\begin{equation}
\transpose{U}\fA = \left[ \begin{array}{cccc} \transpose{u_{1}}\fA & \transpose{u_{2}}\fA & \cdots & \transpose{u_{\nrank}}\fA \end{array} \right] \text{,}
\label{eq_UfA}
\end{equation}
\[ \text{such that \ } \left[\transpose{U}\fA\right](h,\psi) = \left[\transpose{u_{h}}\fA\right](\psi) \text{.} \]
We use the same convention for principal volumes:
\begin{equation}
\transpose{U}\fF = \left[ \begin{array}{cccc} \transpose{u_{1}}\fF & \transpose{u_{2}}\fF & \cdots & \transpose{u_{\nrank}}\fF \end{array} \right] \text{,}
\label{eq_UfAv}
\end{equation}
\[ \text{such that \ } \left[\transpose{U}\fF\right](h,\hk) = \left[\transpose{u_{h}}\fF\right](\hk) \text{.} \]
We will also refer to the Bessel and spherical harmonic representations 
of the principal images and principal volumes as $\transpose{U}\bA$ and $\transpose{U}\sF$, respectively.

\subsection{Volume reconstruction}
\label{sec_Volume_Reconstruction}

Most approaches to estimating molecular volumes in cryo-em involve some kind of reconstruction process which is performed multiple times (over the course of iterative refinement).
For standard AM, the inputs to the reconstruction process are the collection of $\nimage$ picked particle images $\fA_{j}(\vk)$ and associated CTF-functions $CTF_{j}(\vk)$, along with estimates of the viewing angles $\tau_{j}^{\estimated}$ and displacements $\vd_{j}^{\estimated}$.
The output is the estimate $\fF^{\estimated}(\vk)$.

In our case, we will make use of both the standard least-squares 
reconstruction, as well as a more computationally efficient 
approximate Fourier inversion which we refer to as 
`quadrature-back-propagation' (section \ref{sec:qbp} below).

\subsection{Least-squares reconstruction}

The standard strategy for reconstruction typically involves solving a 
linear system in a least-squares sense.
In real-space this least-squares system looks like:
\[ 
\left[
\begin{array}{c}
PSD_{1}(\vx) \star \translation(+\vd_{1}) \circ \left\{ S\circ \rotation(\tau_{1})\circ \right\} \\
\vdots \\
PSD_{\nimage}(\vx) \star \translation(+\vd_{\nimage}) \circ \left\{ S\circ \rotation(\tau_{\nimage})\circ \right\} \\
\end{array}
\right]
\cdot
F(\vx)
\approx
\left[
\begin{array}{c}
A_{1}(\vx) \\
\vdots \\
A_{\nimage}(\vx) \\
\end{array}
\right]
\]
whereas in the Fourier domain it is written as:
\[
\left[
\begin{array}{c}
CTF_{1}(\vk) \odot \translation(+\vd_{1},\vk) \odot \left\{ \fSlice\circ \rotation(\tau_{1})\circ \right\} \\
\vdots \\
CTF_{\nimage}(\vk) \odot \translation(+\vd_{\nimage},\vk) \odot \left\{ \fSlice\circ \rotation(\tau_{\nimage})\circ \right\} \\
\end{array}
\right]
\cdot
\fF(\vk)
\approx
\left[
\begin{array}{c}
\fA_{1}(\vk) \\
\vdots \\
\fA_{\nimage}(\vk) \\
\end{array}
\right]
\]
or, more conveniently,
\begin{equation}
\left[
\begin{array}{c}
CTF_{1}(\vk) \odot \left\{ \fS\circ \rotation(\tau_{1})\circ \right\} \\
\vdots \\
CTF_{\nimage}(\vk) \odot \left\{ \fS\circ \rotation(\tau_{\nimage})\circ \right\} \\
\end{array}
\right]
\cdot
\fF(\vk)
\approx
\left[
\begin{array}{c}
\translation(-\vd_{1},\vk) \odot \fA_{1}(\vk) \\
\vdots \\
\translation(-\vd_{\nimage},\vk) \odot \fA_{\nimage}(\vk) \\
\end{array}
\right]
\text{.}
\label{eq_0lsq_orig}
\end{equation}
In this expression, the $\nimage$ viewing angles $\tau_{j}$ and 
displacements $\vd_{j}$ are assumed to be given as input, as are the 
image- or micrograph-specific CTF functions $CTF_{j}(\vk)$ and 
images $\fA_{j}(\vk)$. 
The unknown (to be determined) is the volume $\fF(\vk)$. 

This least-squares problem is standard in maximum-likelihood estimation, and versions of this problem are solved in many software packages, often via conjugate gradient iteration applied to the 
normal equations.

As written, the least-squares problem above has a block-structure: different $k$ values are independent from one another.
Thus, the problem can be solved by determining 
the $k$-shell $\fF(k,\hk)$ from the $k$-rings $\fA_{j}(k,\psi)$ for one
$k$ value at a time.
However, the problem stated in \eqref{eq_0lsq_orig} is slightly misleading, 
as the residuals for each of the $k$ values are not directly comparable
due to the fact that
different $k$-shells have different variances associated with them 
(see section \ref{sec_Variance_Scaling}).
To account for this, we rewrite \eqref{eq_0lsq_orig} as:
\begin{equation}
\left[
\begin{array}{c}
\sqrt{w(k)}\cdot CTF_{1}(\vk) \odot \left\{ \fSlice\circ \rotation(\tau_{1})\circ \right\} \\
\vdots \\
\sqrt{w(k)}\cdot CTF_{\nimage}(\vk) \odot \left\{ \fSlice\circ \rotation(\tau_{\nimage})\circ \right\} \\
\end{array}
\right]
\cdot
\fF(\vk)
\approx
\left[
\begin{array}{c}
\sqrt{w(k)}\cdot \translation(-\vd_{1},\vk) \odot \fA_{1}(\vk) \\
\vdots \\
\sqrt{w(k)}\cdot \translation(-\vd_{\nimage},\vk) \odot \fA_{\nimage}(\vk) \\
\end{array}
\right]
\text{,}
\label{eq_0lsq_full}
\end{equation}
where the function $w(k_{r})$ is the radial quadrature weight.
With this formulation, the residual (calculated using the standard vector 
2-norm) is (up to a quadrature error) equal to an integral representing 
the log-likelihood that the volume $\fF(\vk)$ could have produced the 
observed images $\fA_{j}(\vk)$ (under the assumptions about the 
noise-model described in section \ref{sec_Image_noise_model}).

\subsection{Principal-mode least-squares reconstruction}

The same general methodology above can be applied to the principal 
images directly.
To describe the modified least-squares system, we first define the 
low-rank decomposition of the radial component of the CTF-functions.
Let $CTF(k,j)$ correspond to the radial component of the CTF-function $CTF_{j}(\vk)$.
Using the discretized analog of a functional SVD we represent
$CTF(k,j)$ by
\begin{eqnarray}
CTF(k,j) & = & \sum_{h'} u'_{h'}(k)\cdot\sigma'_{h'}v'_{h'}(j) \text{,}
\label{eq_CTF_USV}
\end{eqnarray}
which we can truncate at rank $\nrank'$ in accordance with the global 
tolerance.
In the numerical experiments we present below, the image pools typically 
span $30$ or so micrographs, usually resulting in $\nrank'$ about $2-3$.

Using this $\nrank'$-rank decomposition of the CTF-array 
(across the image pool), we can approximately represent the noiseless 
principal images as follows:
\begin{eqnarray}
CTF_{j}(\vk)\odot\left\{\fSlice\circ\rotation(\tau_{j})\circ\fF\right\}(\vk)
& \approx &
\sum_{h'=1}^{\nrank'} v'_{h'}(j) \sigma'_{h'} \cdot u'_{h'}(\|\vk\|) \odot\left\{\fSlice\circ\rotation(\tau_{j})\circ\fF \right\}(\vk) \text{.}
\label{eq_a_UCTF_Y}
\end{eqnarray}
After using the radial principal modes to compress both this approximation and the associated component of the right hand side of \eqref{eq_0lsq_full}, we obtain:
\begin{eqnarray}
\sum_{h'=1}^{\nrank'} v'_{h'}(j) \sigma'_{h'} 
\cdot
\transpose{u_{h}} \left[ u'_{h'}(\vk) \odot\left\{\fSlice\circ\rotation(\tau_{j})\circ\fF \right\}(\vk) \right] 
=
\transpose{u_{h}}\left[ \translation(-\vd_{j},\vk)\odot \fA(\vk) \right]
\text{.}
\end{eqnarray}
By defining the CTF-modulated principal volume:
\[ 
\fG(k,\hk ; h',\fF) = \transpose{U}\cdot\left[u'_{h'}(k)\odot\fF(k,\hk)\right] \text{,}
\]
we can now write the compressed least-squared problem as:
\begin{equation}
\sum_{h'=1}^{\nrank'} \sigma'_{h'}
\left[
\begin{array}{c}
v'_{h'}(1) \left\{ \fS\circ \rotation(\tau_{1})\circ \right\} \\
\vdots \\
v'_{h'}(\nimage) \left\{ \fS\circ \rotation(\tau_{\nimage})\circ \right\} \\
\end{array}
\right]
\cdot
\fG(\vk ; h' , \fF)
\approx
\left[
\begin{array}{c}
\transpose{U}\left[\translation(-\vd_{1},\vk) \odot \fA_{1}(\vk)\right] \\
\vdots \\
\transpose{U}\left[\translation(-\vd_{\nimage},\vk) \odot \fA_{\nimage}(\vk)\right] \\
\end{array}
\right]
\text{,}
\label{eq_0lsq_pm}
\end{equation}
where the radial reweighting of $\sqrt{w(k)}$ in \eqref{eq_0lsq_full} is accounted for by our definition of $\transpose{U}$.
Note that, when $\nrank'=\nrank=\rmax$ the equation above is equivalent to \eqref{eq_0lsq_full}.
However, when $\nrank'<\rmax$ and $\nrank<\rmax$ this will correspond to a projected version of \eqref{eq_0lsq_full}.

Note also that, with this representation, we need not actually solve 
for $\fF(\vk)$.
Instead, we can use the principal images (on the right hand side) to solve the above problem (in a least-squares sense) for the collection of CTF-modulated principal volumes $\fG(\vk;h',\fF)$.
Solving \eqref{eq_0lsq_pm} can be less computationally costly than solving \eqref{eq_0lsq_full} when the total number of degrees of freedom within the various $\fG$ is less than the total number of degrees of freedom within $\fF$ (i.e., when $\nrank\times\nrank'<\rmax$).
Even when this is not the case, the CTF-modulated principal volumes $\fG$ (the solution to \eqref{eq_0lsq_pm}) can still be more accurate (in terms of comparison to the ground truth) than the solution $\fF$ to \eqref{eq_0lsq_full}.
We see this improvement most acutely when the principal mode reduction $\transpose{U}$ effectively removes noise from the images $\fA_{j}$ (e.g., if there are many frequency rings which contain little to no information).

We remark that, when dealing with least-squares reconstruction and principal mode reduction, `the diagram commutes'. 
That is, because the principal mode projection $\transpose{U}$ is orthonormal, the principal mode projection $\transpose{U}\left\{CTF\cdot\fF\right\}$ (i.e., the principal projection of the full volume which solves the original least-squares problem in \eqref{eq_0lsq_full}) will be the same as the solution $\sum_{h'}v'\sigma'\fG$ to the principally-projected least-squares problem in \eqref{eq_0lsq_pm}, provided of course that the $\nrank'$-rank decomposition of the CTF-array is accurate (e.g., if $\nrank'$ is taken to $\rmax$). 
Thus, if $\nrank\times\nrank'$ is indeed larger than $\rmax$, then instead of solving \eqref{eq_0lsq_pm}, we can obtain $\fG$ more easily by just solving the original least-squares problem for $\fF$ and then projecting down to the CTF-modulated principal volumes $\fG$.

Note also that, due to linearity, the noiseless images associated with the CTF-modulated principal volumes $\fG$ will approximately equal the noiseless principal images associated with the volume $\fF$, with errors controlled by the truncation $\nrank'$ of the CTF-values.
That is to say, a principal mode compression of \eqref{eq_a_UCTF_Y} can be written as: 
\begin{eqnarray}
\transpose{U}\left[ CTF_{j}(\vk)\odot\left\{\fS\circ\rotation(\tau_{j})\circ\fF\right\}(\vk) \right]
& \approx &
\sum_{h'=1}^{\nrank'} v'_{h'}(j) \sigma'_{h'} \cdot \left\{ \fS\circ\rotation(\tau_{j})\circ\fG(\vk;h',\fF) \right\}(\vk) \text{.}
\label{eq_a_UCTF_UX_Y}
\end{eqnarray}
This relationship can be used to align the principal images $\transpose{U}\fA_{j}$ directly to the solution $\fG$, without necessarily needing to reconstruct $\fF$ (see section \ref{sec_Image_Alignment}).

\paragraph{Reconstruction via quadrature back propagation (QBP)}
\label{sec:qbp}

In our experience, an accurate approximation of the solution to the least-squares problem above can require a significant amount of computation time when the conjugate gradient method requires many iterations to converge.
As a more efficient (but typically less accurate) alternative, we employ a slightly different strategy.
We refer to this alternative as `quadrature-back-propagation' (QBP).

The basic observation underlying QBP is that, assuming that there was no image noise and no CTF-attenuation, then the volume and images would solve:
\[
\left[
\begin{array}{c}
\fSlice\circ \rotation(\tau_{1})\circ  \\
\vdots \\
\fSlice\circ \rotation(\tau_{\nimage})\circ  \\
\end{array}
\right]
\cdot
\fF
=
\left[
\begin{array}{c}
\translation(-\vd_{1},\vk)\odot\fA_{1}(\vk) \\
\vdots \\
\translation(-\vd_{\nimage},\vk)\odot\fA_{\nimage}(\vk) \\
\end{array}
\right]
\]
with each template $\fSlice\circ\rotation(\tau_{j})\circ\fF$ corresponding to a cross-section of the volume $\fF(\vk)$.
Each of these cross-sections can be viewed as a sample of the overall volume $\fF(\vk)$, with values given by the appropriate entry on the right-hand-side of the equation above.
With sufficiently many cross-sectional samples, a functional representation of the volume $\fF$ can be reconstructed via quadrature in the spherical harmonic basis via a spherical harmonic transform. That is, we have the approximate
formula
\[ \sF(k;l,m) \approx \sum_{p} w_{p}\left[Y_{l}^{m}(\hk_{p})\right]^{\dagger}\fF(k,\hk_{p})\ \ \forall \ \ l,\|m\|,l',\|m'\| \leq \lmax\text{.} \]
where $\{ \hk_p,w_p \}$ correspond to a quadrature rule on the sphere.
The number of such nodes is of the order 
$[\bigO(k)]^2$ on the spherical $k$-shell.
(Equispaced points in the azimuthal direction and Legendre points in 
the polar angle, for example, are spectrally accurate.)
Because the data points will not generally coincide with the necessary
quadrature nodes, we interpolate from the given data to the nodes 
in some local neighborhood.
A simple rule would be to take the local average over all data points
that lie within some distance of the node.
Note that this strategy is not iterative and is 
typically an order of magnitude (or so) faster than a least-squares 
reconstruction.

In the absence of image-noise, the error in the reconstruction of $\sF(k;l,m)$ will come from the `interpolation error' associated with interpolating from the cross-sectional sample data to the quadrature nodes. 
When the image-noise is nonzero and the interpolation regions 
are fixed in size, then the approximation to $\fF(k,\hk_{p})$ will exhibit two sources of error.
In addition to the interpolation error, there is a sampling error
due to the image noise at each cross-sectional sample-point.
In practice, we generally choosing the local neighborhood of the quadrature 
nodes to be sufficiently large to contain roughly $\bigO(10^{1-2})$ sample points. This is usually sufficient for two digits of accuracy (coinciding with our usual global error tolerance).

The above strategy can easily be generalized to deal with image-specific CTF-values.
The QBP approach above essentially treats each sample within the local
interpolation region of each quadrature node as a noisy observation of 
$\fF(k,\hk_{p})$.
The average can be thought of as a local maximum-likelihood estimate.
Since different data points might correspond to different CTF functions,
we can simply compute the CTF-weighted local averages instead.
This strategy coincides with the CTF-weighting of Bayesian inference 
described in \cite{Scheres2012,Scheres2012b}.

\subsection{Remark regarding further template compression}

If the ranks $\nrank$ and $\nrank'$ are both sufficiently small 
(e.g., $2-4$ and $1-2$, respectively), then the dimensionality of the 
space of CTF-modulated principal templates (as measured via the numerical rank) is very small, and usually quite a bit smaller than the actual number of principal templates (i.e., the number of viewing-directions). 
This would allow for additional compression, using a reduced set 
of `basis'-templates (which span the space of templates) rather than the actual templates themselves. 
In our numerical experiments we do not implement this compression because our $\nrank$ and $\nrank'$ are typically too large for it to be efficient.

\subsection{Image alignment}
\label{sec_Image_Alignment}

Just as for the reconstruction step, image alignment is often performed 
multiple times during molecular refinement.
For standard AM, the inputs to the alignment process are the estimated 
molecular volume $\fF^{\estimated}(\vk)$, as well as the picked-particle images $\fA_{j}(\vk)$ and associated CTF-functions $CTF_{j}(\vk)$. 
The outputs are the estimated viewing angles $\tau_{j}^{\estimated}$ and displacements $\vk_{j}^{\estimated}$ for the various images.

In our case we will make use of both the standard maximum-likelihood alignment, as well as a more stable `maximum-entropy' alignment (described below).

\subsection{Alignment via maximum likelihood}

For this, we calculate the quantities $\cX$, $\cY$ and $\cZ$, 
described in section \ref{sec:correlations},
for a collection of $\tmax=\bigO(\kmax^{2})$ viewing directions $\eazimub_{t}$ and $\epolara_{t}$ corresponding to our spherical shell quadrature grid 
for $\hk$ at $k=k_{\rmax}\approx\kmax$.
For each of these viewing-directions $(\eazimub_{t},\epolara_{t})$, we consider a range of $O(\kmax)$ values for $\egammaz$ in $[0,2\pi)$.
For the displacements, we typically limit them to lie in a small disk of 
radius $\vdlocal\sim 0.01$ or so, taking advantage of the low numerical rank of the set of translation-operators within this disk 
\cite{RSAB20,Rangan22}.

Given the correlations $\cZ(\eazimub_{t},\epolara_{t};\fA_{j})$, a common strategy for updating the viewing angles $\tau_{j}^{\estimated}$ and displacements $\vd_{j}^{\estimated}$ (referred to as `maximum-likelihood' alignment), involves selecting the $\tau_{j}^{\estimated}$ (and, as a result, the $\vd_{j}^{\estimated}$) for each image $j$ to maximize $\cZ(\eazimub_{t},\epolara_{t})$ over the $\tmax$ viewing directions for that particular image $j$.
In the ideal situation where the image formation model is accurate, this choice corresponds to the viewing-angle that is most likely to have produced the particular image $\fA_{j}$, given that the molecule is indeed $\fF^{\estimated}$.
In practice we do not define $\tau_{j}^{\estimated}$ to be the very best $\eazimub$ and $\epolara$ for each image, instead randomly selecting an $\eazimub_{t}$ and $\epolara_{t}$ for which $\cZ(\eazimub_{t},\epolara_{t})$ is within the top $5^{\text{th}}$-percentile of its values for that image \cite{Barnett2017}.
The results of our numerical-experiments are quite insensitive to this value of $5\%$; any choice from $1\%$ to $\sim 10\%$ achieves similar performance for maximum-likelihood alternating minimization.

Let us define the function $\sort[\cdot]$ which maps lists of real-numbers to lists of positive integers such that, for any vector 
$\vect{a}\in\Real^{L}$, the value of $\sort[\vect{a}]_{l}$ is 
the (integer) rank of $a_{l}$ (ranging from $1,\ldots,L$ as it stands within the entries of $\vect{a}$).
As noted in section \ref{sec:correlations}, 
we can see that the template ranks are obtained by applying $\sort[\cdot]$ to each `column' of the array $\cZ(\eazimub_{t},\epolara_{t};\fA_{j})$ (fixing $j$), while the image ranks are obtained by applying $\sort[\cdot]$ to each `row' of the same array (fixing $t$).
The maximum-likelihood alignment described above depends on the template ranks, while the maximum-entropy alignment (described below) depends
on the image ranks.

\subsection{Alignment via maximum-entropy}

The maximum-likelihood alignment described above can be thought of as a `greedy' algorithm which selects the estimated viewing angles for each image without consideration of the choices made for the other images.
As a result of this greedy strategy, the estimated viewing angle distribution (across all $\tau_{j}^{\estimated}$) is not controlled, and can end up becoming quite far from the true viewing-angle distribution (until the 
volume used to generate the templates is sufficiently accurate).

To avoid this phenomenon, one can choose the $\tau_{j}^{\estimated}$ while constraining the distribution of estimated viewing angles in some way.
A simple special case of this strategy is to constrain the distribution of estimated viewing angles to be uniform.
We refer to such an alignment strategy as `maximum-entropy' alignment 
(referring to the entropy of $\mu_{\tau}^{\estimated}$).
In practice we implement this maximum-entropy alignment as follows:
\begin{enumerate}
    \item Declare all the images `unassigned'.
    \item Define a list of $\nimage$ uniformly distributed viewing-directions $\hk_{j}=\left(\eazimub_{j},\epolara_{j}\right)$ on the sphere.
    (This is typically done by cycling through a randomly-permuted list of the $\tmax$ viewing-angles $(\eazimub_{t},\epolara_{t})$ until $\nimage$ viewing-directions have been selected).
    \item Step through the list of viewing directions.
    For each viewing-direction (defined by a particular $\eazimub$ and $\epolara$), assign the image $\fA_{j}$ that has not previously been assigned and which maximizes $\cZ(\eazimub,\epolara;\fA_{j})$ (over the unused image-indices $j$).
    When assigning an image $\fA_{j}$ to a particular viewing direction $\eazimub$, $\epolara$, we set the in-plane rotation $\egammaz_{j}^{\estimated}$ and displacement $\vd_{j}^{\estimated}$ implicitly via the $\argmax$ (see the definition of $\cZ$ in Eq \ref{eq_cZ}).
\end{enumerate}

Note that the third step of this procedure uses the image ranks for each template (rather than the template ranks for each image).
In some cases the image ranks and template ranks contain similar information (see, e.g., the results for ISWINCP in Fig \ref{Fig_viewing_angle_distribution_SET_ONE}).
However, it is more typical for the image ranks to be more informative than the template ranks 
(see the other case studies in Figs. \ref{Fig_viewing_angle_distribution_SET_ONE} and \ref{Fig_viewing_angle_distribution_SET_TWO}).

Intuition might suggest that maximum-entropy alignment is a reasonable 
strategy if the true viewing angle distribution is known to be uniform,
and not otherwise.
Nevertheless, as illustrated in the main text, our EMPM iteration still often produces results that are significantly better than maximum-likelihood alignment, even for cryo-EM data where the true viewing angle distribution is rather nonuniform.
We study this behavior in some detail for simplified models 
in sections \ref{sec_MRA} and \ref{sec_MSA} below.

\section{Multi-Reference Alignment (MRA)}
\label{sec_MRA}
Multi-Reference Alignment (MRA) is an abstraction of
the cryo-EM molecular reconstruction problem, reduced to a two-dimensional
setting.

In this context the image displacements and CTF-modulation are typically ignored, and we assume (i) the imaged molecule itself is translationally invariant in the $x_{3}$-direction, and
(ii) the imaged molecule is only observed from the $x_{3}$-axis.
Because of (i), the function $\xF(\vx)$ will be determined by $\xF(x_{1},x_{2})$, and $\fF(\vk)$ will be restricted to the equatorial-plane, as denoted by $\fF(k,\psi)$.
Additionally, due to (ii), each `image' $\fA_{j}$ will be a noisy version of $\fF(k,\psi)$, subject to an unknown in-plane rotation $\egammaz_{j}$:
\[ \fA_{j}(k,\psi) = [\rotation(\gamma_{j})\circ\fF](k,\psi) + \text{noise}\ = \fF(k,\psi-\gamma_{j}) + \text{noise}.  \]
If we further assume that the true volume $\fF$ involves only a 
single $k$-shell, then each image is restricted to a single $k$-ring, 
and can be considered a periodic function on $[0,2\pi]$.

\begin{figure}
    \centering
    \includegraphics[width=5.8in]{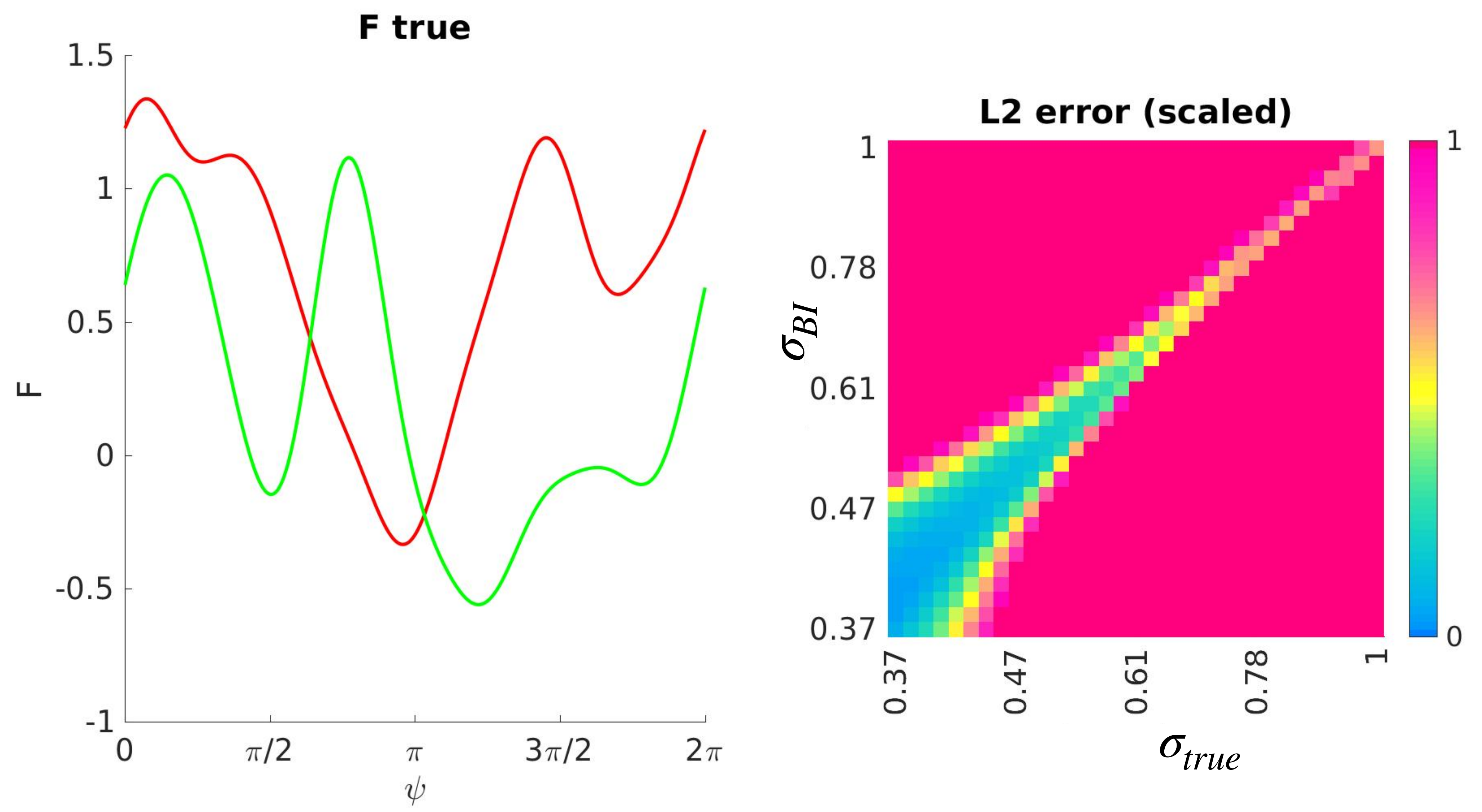}
    \caption{
In this figure we illustrate the typical errors observed when using 
Bayesian inference to tackle the MRA problem.
For this numerical experiment we set 
$\fG^{\true}$ to be a random function with $\qmax=8$ 
(i.e., $1+2\qmax=17$ degrees of freedom).
The real (red) and imaginary (blue) parts of $\fG^{\true}$ are shown on 
the left. We set $\nimage=2048$, with $\psi_{j}^{\true}$ uniformly 
distributed in $[0,2\pi]$.
    We vary the true noise-level (horizontal) and the estimated noise-level (vertical).
    On the right-hand-side we show the $L^2$-error between $\fG^{\true}$ and $\fG^{\estimated}$.
    We calculate $\fG^{\estimated}$ via bayesian-inference, starting with a random function with the same number of degrees of freedom as $\fG^{\true}$, and iterating until convergence.
    The $L^2$-error on the right is scaled so that `$1$' corresponds to the error between $\fG^{\true}$ and a constant function with the same mean.
    Note that when the true noise-level is moderate (e.g., around $0.5$) there is only a narrow range of estimated noise-levels which admit accurate reconstructions (i.e., when $\sigma^{\estimated}$ is close to $\sigma^{\true}$).
    Outside this narrow band the reconstruction can be quite poor, even when $\sigma^{\estimated}$ is only $20\%$ higher or lower than $\sigma^{\true}$.
    }
    \label{Fig_MRA_sheres}
\end{figure}

The MRA problem has been well studied 
\cite{bendory_mra_2021,bendory_mra_2022,perry_mra_2019}, 
and can be used to gain insight into the behavior of many 
reconstruction strategies. It can be used, for example, to investigate 
the sensitivity of Bayesian inference (BI) to mismatches between
the estimated and true noise levels (see Fig \ref{Fig_MRA_sheres}).

\section{Multi-Slice Alignment (MSA)}
\label{sec_MSA}

An `orthogonal' sub-problem to MRA, which we refer to as Multi-Slice 
Alignment (MSA), can be described as follows.
As in MRA, we ignore displacements and the CTF and assume that the volume is invariant in the $x_{3}$-direction.
However, in MSA, we assume that
the imaged molecule is only observed from the equatorial plane.
With this assumption the images will correspond to projections of $\xF(\vx)$ onto one-dimensional
lines in the $x_{1}$-$x_{2}$-plane (i.e., samples of the two-dimensional 
Radon transform of $\xF(x_{1},x_{2})$).
Thus, retaining the notation $\fA_{j}$, each image corresponds to 
a one-dimensional slice of $\fF(k,\psi)$, expressed as 
a function of $k\in\Real$ for a specific (but unknown) $\psi_{j}$:
\[ \fA_{j}(k) = \fF(|k|,\psi_{j}) + \text{noise} \text{\ if\ } k>0 \text{,\ or\ \ } \fF(|k|,\psi_{j}+\pi) + \text{noise} \text{\ if\ } k<0 \text{.} \]
To simplify the scenario further, we restrict each image to the ray of positive $k$:
\[ \fA_{j}(k) = \fF(k,\psi_{j}) + \text{noise} \quad\text{for $k>0$} \text{.} \]

\begin{figure}
    \centering
    \includegraphics[width=5.5in]{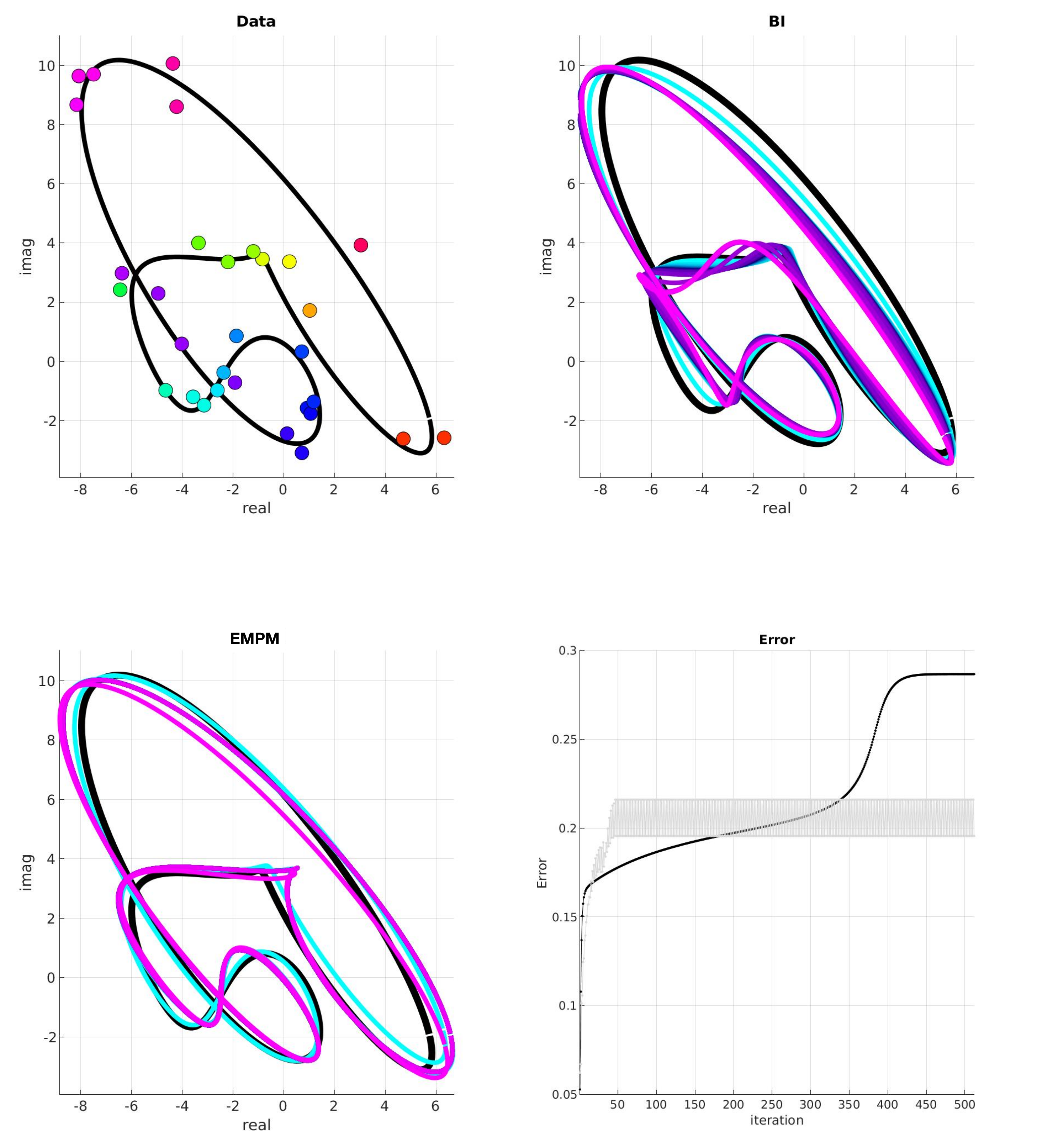}
    \caption{
    In this figure we illustrate the MSA problem with a single
shell at $k=1$.
    In the first subplot we show one example of a function $\fG^{\true}(\psi)$ plotted as a closed curve in $\Complex$ (black).
    $\nimage=32$ images are sampled from $\fG(\psi)$, with each $\psi_{j}^{\true}$ drawn independently from a uniform distribution on $[0,2\pi]$ (circles, colored periodically according to $\psi_{j}^{\true}$).
    The true noise-level used when generating the images is $\sigma^{\true}=0.5$.
    Starting with the ground-truth $\fG^{\initial}=\fG^{\true}$, we apply Bayesian inference (BI) with the estimated noise level $\sigma^{\BayesianInference}$ set to equal $\sigma^{\true}$.
    The resulting $\fG^{\BayesianInference}$ are shown in the second subplot, with early iterations in cyan and later iterations in magenta.
    We also apply the EMPM strategy, again starting with the ground-truth; 
the resulting $\fG^{\EMPM}$ are shown in the third subplot.
    Errors for BI (black) and EMPM (gray) are shown in the fourth sub-plot.
    Note that the error for EMPM oscillates as the iterations themselves alternate between maximum-likelihood- and maximum-entropy-alignment.
    In Fig \ref{Fig_MSA_sheres_lambda} below we use as the error for EMPM 
the higher of the two alternating values.
    The data here correspond to one trial with $\sigma=0.5$ and 
$\lambda=1$ in Fig \ref{Fig_MSA_sheres_lambda}.
    }
    \label{Fig_MSA_sheres_illustration}
\end{figure}

If we further assume that the true volume $\xF$ can take on arbitrary complex values, and that $\fF$ involves only a single $k$-shell, then $\fF$ corresponds to a closed-curve in $\Complex$, written as the function $\fG(\psi)$, parametrized by $\psi$.
Each image is now simply a point, corresponding to a noisy observation of $\fG(\psi)$.
If we assume that the true volume $\fF$ involves a discrete set of $\rmax$ $k$-values $\{k_{1},\ldots,k_{\rmax}\}$, then $\fG$ will be a closed curve in $\Complex^{\rmax}$ (once again parametrized by $\psi$), and each image will be a vector in $\Complex^{\rmax}$ corresponding to a noisy observation of $\fG(\psi_{j})$.
An example illustrating the MSA problem for $\rmax=1$ is shown in 
Fig \ref{Fig_MSA_sheres_illustration}. 

\paragraph{Numerical experiments:}

We believe that MSA can be a useful testbed for studying the 
performance of various strategies for cryo-EM reconstruction.
As a simple example, let $\rmax\equiv 1$ above and define 
$\fG^{\true}:[0,2\pi)\rightarrow\Complex$ using a fixed number of 
Fourier modes:
\[ \fG^{\true}(\psi) = \sum_{q=-\qmax}^{+\qmax} \bG^{\true}(q)\cdot\exp(iq\psi) \text{,} \]
with the shape of $\fG^{\true}$ determined by the $2\qmax+1$ 
coefficients $\bG^{\true}_{q}$.
In the experiments below, we generate $\fG^{\true}$ randomly by drawing 
each of the $\bG^{\true}_{q}$ from a (complex) Gaussian distribution with 
variance $1$.

Having drawn a particular $\fG^{\true}$, we sample this 
molecule at $\nimage$ points, producing a collection of 
simulated data points $\fA_{j}$ defined by:
\[ \fA_{j} = \fG^{\true}(\psi_{j}) + \epsilon_{j}\text{,} \]
where the $\psi_{j}$ are drawn independently and uniformly from $[0,2\pi)$, and each $\epsilon_{j}$ is an independent complex random variable drawn from a normal distribution with variance $\sigma$.

Now, given the collection of $\nimage$ images $\{\fA_{j}\}$, we can try to reconstruct $\fG^{\true}$ using the strategies discussed in the main text.
Due to the simplicity of MSA, much of the apparatus described in the main text simplifies dramatically. For example, reconstruction involves a subsampled  and/or non-uniform Fourier transform, and alignment involves calculating the distance between each image point and the curve $\fG$ as a function of $\psi$.

Many features of  this toy problem are expected.
For example, when $\sigma=0$ the images lie exactly on the 
curve $\fG^{\true}$.
If we begin the alignment process with an initial guess 
$\fG^{\initial}=\fG^{\true}$ for the molecule, then each image will be 
aligned exactly to its correct viewing-angle $\psi_{j}^{\true}$ 
(up to a global phase-shift).
Subsequent reconstruction steps will 
produce $\fG^{\estimated}=\fG^{\true}$, so long as the number of 
observed points $\nimage$ is at least as large as the number of 
modes $1+2\qmax$, and the reconstruction operator is well-conditioned 
(i.e., if the image-points are not too closely clumped together).
Put more simply, if $\sigma=0$ and $\nimage\geq 1+2\qmax$ then, 
generically speaking, a maximum-likelihood alternating minimization 
starting with $\fG^{\initial}=\fG^{\true}$ as an initial-function will 
be stable.

However, even in this simple scenario (i.e., with $\sigma=0$), 
if we were to attempt maximum-likelihood alternating minimization 
with a different initial function then we are not guaranteed to 
converge to $\fG$.
Indeed, if $\nimage$ is not too large, there will typically be many 
different functions with the same bandlimit of $\qmax$ that pass close 
to the observed points $\fA_{j}$.
Which of these functions we converge to will depend on the details of 
the reconstruction process.
The situation becomes more complicated once we add image noise; 
when $\sigma>0$ even an initial guess of $\fG^{\initial}=\fG^{\true}$ 
will not necessarily lead to convergence to a function $\fG^{\estimated}$ 
that is close to $\fG^{\true}$.  

To probe this phenomenon, we perform a set of numerical-experiments.
We initialize our reconstruction with a function of the form:
\[ \fG^{\initial} := \lambda\cdot \fG^{\true} + (1-\lambda)\cdot \fG^{\random} \text{,} \]
where $\fG^{\random}$ is a random function drawn the same way as (but independent from) $\fG^{\true}$, and $\lambda$ is a parameter describing how close $\fG^{\initial}$ is to the ground truth.
Given a particular $\fG^{\initial}$, we perform the reconstruction, iterating until the resulting $\fG^{\estimated}$ converges.
For each trial we can then measure the error $E$ between $\fG^{\estimated}$ and $\fG^{\true}$.

In order to reward trials where $\fG^{\estimated}$ and $\fG^{\true}$ are close to one another, even if they are parametrized differently, we define the error $E$ as follows:
First we define $d(\fG^{\estimated}\leftarrow\fG^{\true}(\psi))$ to be the distance between $\fG^{\true}(\psi)$ (i.e., a point in $\Complex$) and the curve $\fG^{\estimated}$.
We then define the error $E(\fG^{\estimated}\leftarrow\fG^{\true})$ as the integral of $d^{2}(\fG^{\estimated}\leftarrow\fG^{\true}(\psi))$ over $\psi$.
Similarly, we define $d(\fG^{\true}\leftarrow\fG^{\estimated}(\psi))$ to be the distance between $\fG^{\estimated}(\psi)$ and the curve $\fG^{\true}$, and define the error $E(\fG^{\true}\leftarrow\fG^{\estimated})$ as the integral of $d^{2}(\fG^{\true}\leftarrow\fG^{\estimated})$ over $\psi$.
Finally, we define the overall error $E=\max(E(\fG^{\true}\leftarrow\fG^{\estimated}),E(\fG^{\estimated}\leftarrow\fG^{\true}))$ to be the larger of the two.

\begin{figure}
    \centering
    \includegraphics[width=5.5in]{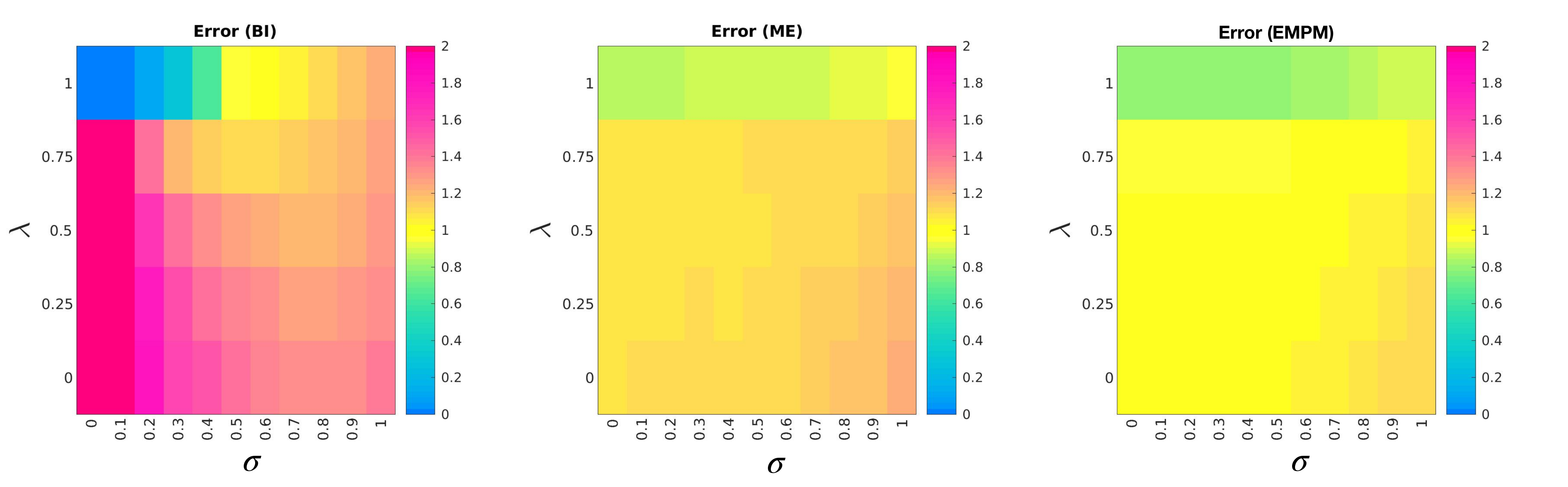}
    \caption{
    In this figure we illustrate the recovery of Bayesian inference (left), 
maximum-entropy alternating minimization (center) and EMPM
as applied to the MSA problem (see text).
For each trial in this numerical experiment we define the true molecule 
using a $\fG$ constructed with $\qmax=8$, corresponding to $1+2\qmax=17$ 
degrees of freedom (each of the coefficients of $\fG$ are drawn 
independently from a normal distribution).
We sample $\nimage=32$ images from $\fG$ (i.e., roughly $2$ times the 
number of modes) and add noise with a magnitude given by $\sigma$ 
ranging from $0$ to $1$ (defined relative to the norm of $\fG$) 
(horizontal).
For each noise level $\sigma$ we perform each reconstruction starting 
with an initial function of the form 
$\lambda\cdot\fG + (1-\lambda)\cdot\fG^{\noise}$, 
with $\fG^{\noise}$ drawn randomly (and independently from $\fG$) 
and $\lambda$ ranging between $0$ and $1$ (vertical).
For each value of $(\sigma,\lambda)$, we simulate $512$ randomly 
generated trials; the heatmaps above average over these trials.
Bayesian inference (left) is run 
with estimated noise level equal to the true noise level.
In the middle, we show results for a simple alternating-minimization 
using only maximum-entropy alignment, rather than maximum-likelihood.
(Maximum-likelihood-based AM fared even more poorly.)
On the right we show the results for EMPM,
with the error measured at the final maximum-entropy alignment 
step (which is typically higher than the error after the last
maximum-likelihood alignment step).
    }
    \label{Fig_MSA_sheres_lambda}
\end{figure}

Figure \ref{Fig_MSA_sheres_lambda} illustrates the average $E$ (accumulated over $512$ random trials) for three different reconstruction strategies across a range of $\lambda$ and $\sigma$ values.
On the left we show the results for Bayesian inference, 
where we set the estimated noise level equal to the true noise level 
(i.e., $\sigma$) when reconstructing the molecule.
Note that Bayesian inference is quite accurate when the initial guess is accurate (i.e., when $\lambda=1$).
However, when $\lambda<1$ and there is a mismatch between the initial guess and the true function, Bayesian inference typically performs quite poorly.

Similarly, we can apply standard maximum-likelihood 
alternating minimization (AM).
This is equivalent to Bayesian inference with an estimated noise level 
of $0$.
The results for AM are identical to those shown on the left when $\sigma=0$.
However, as soon as a little image noise is added ($\sigma>0$), 
AM typically converges to a $\fG^{\estimated}$ which is wildly different from $\fG^{\true}$, with errors that are an order of magnitude larger than those shown for Bayesian inference (which was run with an estimated noise-level of $\sigma$).
We do not display the (large) errors observed for standard AM.

In the middle, we show the results for a simple alternating minimization 
scheme using maximum-entropy alignment at each step.
On the right we show the results for our preferred EMPM strategy: namely, alternating between maximum-likelihood  and maximum-entropy alignment.
Note that, while our preferred strategy is less accurate than Bayesian 
inference when $\sigma$ is small and $\fG^{\initial}=\fG^{\true}$, it is more accurate in this regime when $\sigma>0$ and $\lambda<1$.

\begin{figure}
    \centering
    \includegraphics[width=5.5in]{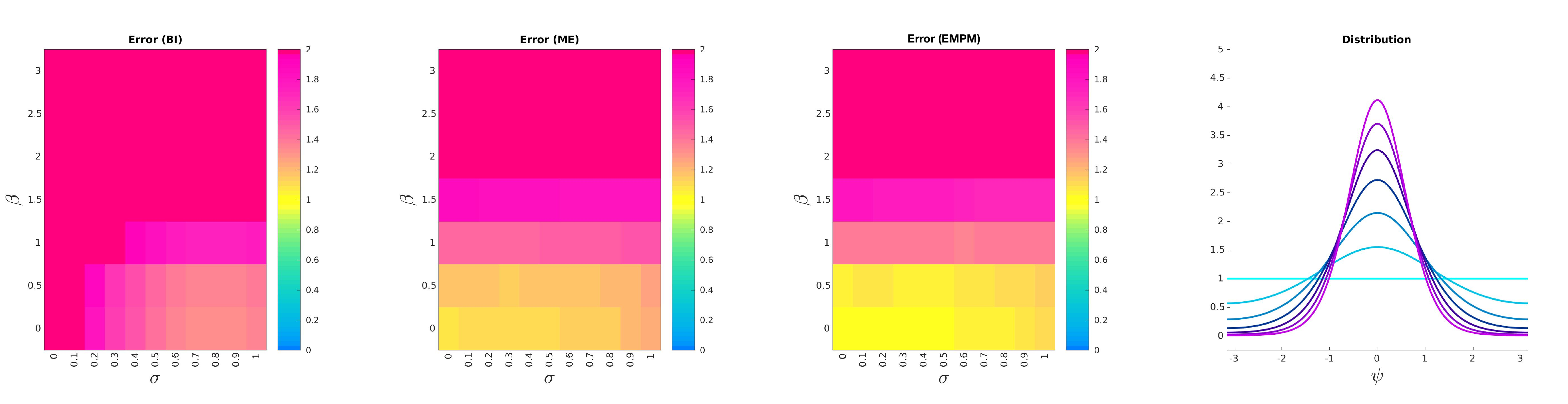}
    \caption{
As in Fig \ref{Fig_MSA_sheres_lambda}, we compare various reconstruction
methods for the MSA problem, but here we focus on viewing angle 
distribution.
For this, we introduce a new parameter $\beta\in[0,3]$, 
where $\beta$ controls the nonuniformity of the viewing angle,
with samples $\psi_{j}$ drawn from the distribution: 
$\exp(\beta\cdot\cos(\psi))/(2\pi{\cal I}_{0}(\beta))$.
Examples of this distribution (multiplied by $2\pi$) are shown on the 
far right.
When $\beta=0$ the viewing angles are uniformly distributed (cyan).
When $\beta=1$ the viewing angle distribution ranges in density 
(from highest to lowest) by a factor of $e^{2}$.
When $\beta=3$ the viewing angle distribution ranges in density by 
a factor of $e^{6}$, and is mostly supported in an interval with 
diameter $\sim\pi$ (magenta).
We fix $\lambda\equiv 0$ (i.e., we start with a random function), and illustrate the errors using the same scale as in Fig \ref{Fig_MSA_sheres_lambda}, this time varying $\beta$ along the vertical axis.
The horizontal axis refers to the noise level as in 
Fig. \ref{Fig_MSA_sheres_lambda}.
    }
    \label{Fig_MSA_sheres_beta}
\end{figure}

Figure \ref{Fig_MSA_sheres_beta} illustrates the results
of a similar numerical experiment, where the images are drawn from a 
non-uniform distribution of viewing angles.
For this numerical experiment, we consider viewing angle distributions of 
the form $\exp(\beta\cos(\psi))/(2\pi{\cal I}_{0})$, 
controlled by a parameter $\beta$.
When $\beta=0$ the viewing angle distribution is uniform, and as $\beta$ increases, the viewing angle distribution becomes more and more sharply peaked around $\psi=0$.
Results for $\lambda=0$ are shown in Fig \ref{Fig_MSA_sheres_beta}.
The results with $\beta=0$ are the same (in distribution) as the $\lambda=0$ scenario in Fig \ref{Fig_MSA_sheres_lambda}.
Note that EMPM is typically more accurate than Bayesian inference in this 
regime.

\section{Relion De-Novo parameters:}
\label{sec_RelionDeNovo}

When using Relion's de-novo molecular reconstruction, we use the following parameters.

\vspace{.2in}

\begin{table}[h]
\begin{center}
\begin{tabular}{|l|l|}
\hline
 {\tt \small --sgd\_ini\_iter 50 } & {\tt \small --sgd\_inbetween\_iter 200 } \\ 
{\tt \small --sgd\_fin\_iter 50 } &  {\tt \small --sgd\_write\_iter 20 } \\
{\tt \small --sgd\_ini\_resol 35 } &  {\tt \small --sgd\_fin\_resol 15 } \\
{\tt \small --sgd\_ini\_subset 100 } &  {\tt \small --sgd\_fin\_subset 500 } \\
 {\tt \small --sgd } &  {\tt \small --denovo\_3dref } \\ 
 {\tt \small --ctf } &  {\tt \small --K 1 } \\
 {\tt \small --sym C1 } &  {\tt \small --flatten\_solvent } \\ 
 {\tt \small --zero\_mask } &  {\tt \small --dont\_combine\_weights\_via\_disc } \\
 {\tt \small --pool 3 } &  {\tt \small --pad 1 } \\  
{\tt \small --skip\_gridding } &  {\tt \small --oversampling 1 } \\
 {\tt \small --healpix\_order 1 } &  {\tt \small --offset\_range 8 } \\ 
 {\tt \small --offset\_step 2 } &  {\tt \small --j 6 } \\
 \hline
\end{tabular}
\end{center}
\end{table}

\vspace{.2in}

Additionally, we set the {\tt \small --particle\_diameter } to be commensurate with the true molecule diameter (rounded up to the nearest $50$ or $100$\AA).

\begin{figure}
    \centering
    \includegraphics[width=5.5in]{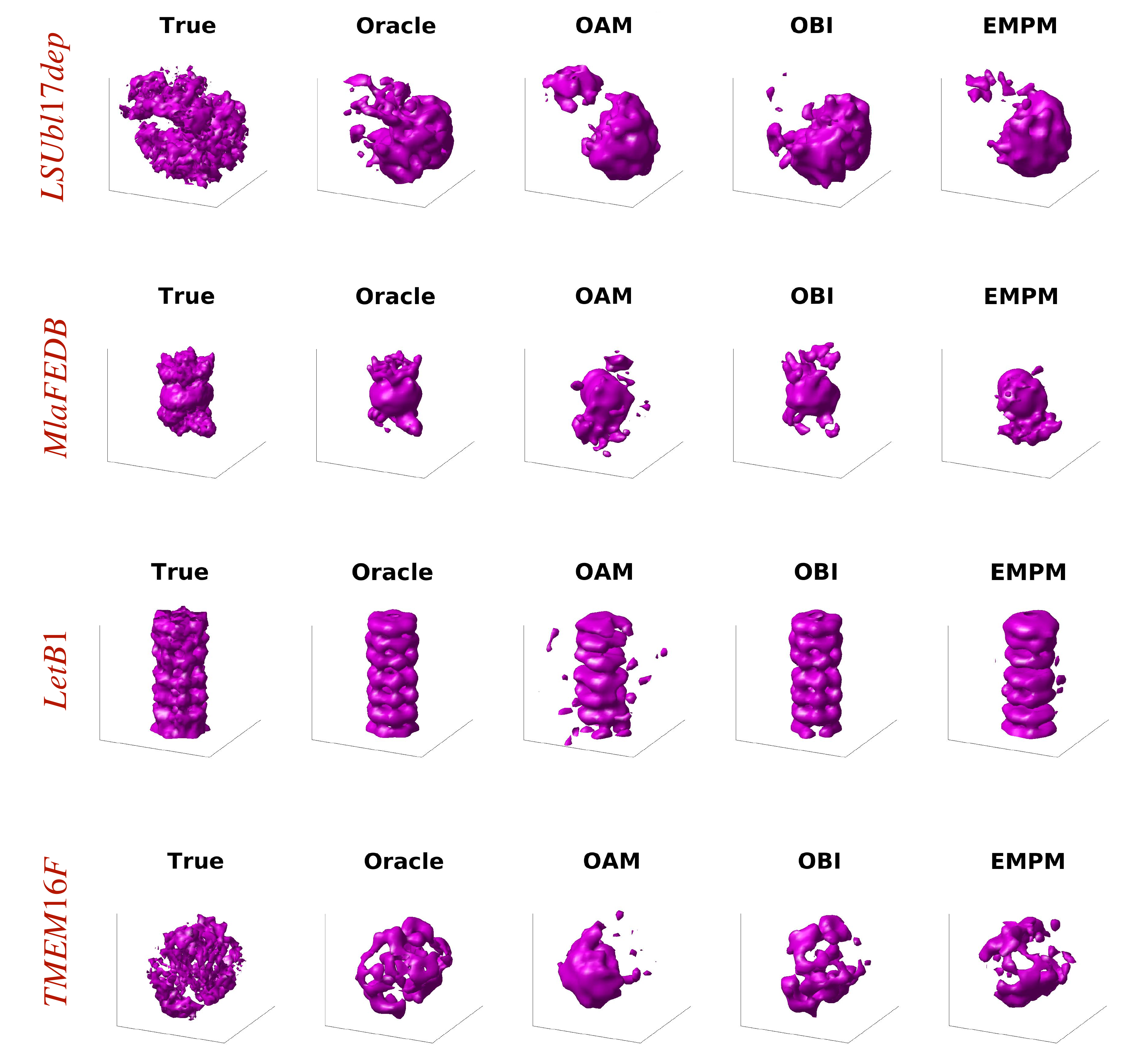}
    \caption{
{
    This figure (analogous to Fig. \ref{Fig_vol_SET_ONE}) plots level sets
from the reconstructions using various 
    Level sets from the volumes $\fF^{\true}$, $\fF^{\Oracle}$, $\fF^{\OracleAMStable}$, $\fF^{\OracleBayesianInference}$, $\fF^{\OracleBayesianInferenceStable}$, and $\fF^{\EMPM}$ (left to right) for
    {\em LSUbl16dep, MlaFEDB, LetB1} and {\em TMEM16F}.
    }}
    \label{Fig_vol_SET_TWO}
\end{figure}

\begin{figure}
    \centering
    \includegraphics[width=5.5in]{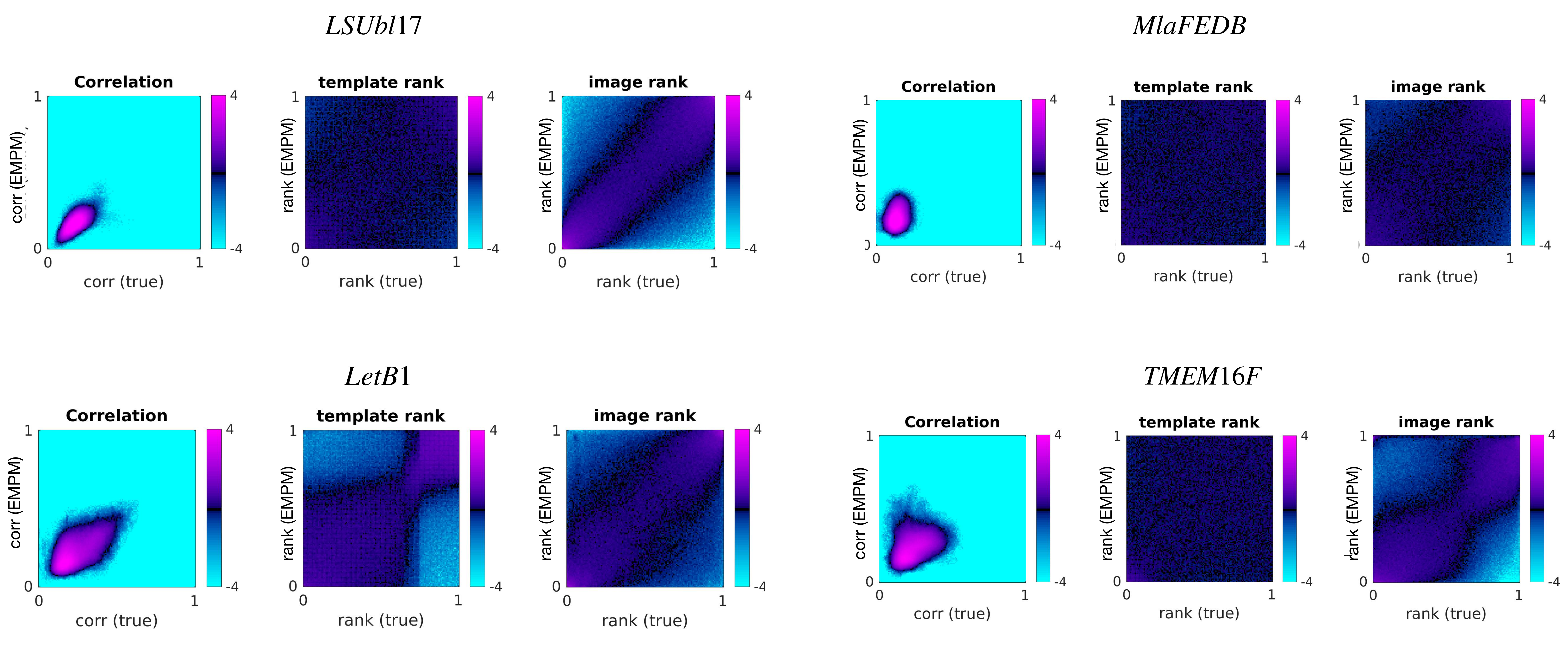}
    \caption{
    The correspondence between estimated and ground-truth viewing angles 
for the reconstructions shown in Fig \ref{Fig_vol_SET_TWO}.
Three panels are shown for each molecule of interest. 
The left-most subplot
is a scatterplot of the image-template correlations. 
The middle subplot is a scatterplot of the template
ranks aggregated over all images, and
the right-most 
subplot is a scatterplot of the image ranks aggregated over all templates.
All are presented as heatmaps, with color corresponding to the
log-density of the joint-distribution.
This figure is the analog of Fig. 
\ref{Fig_viewing_angle_distribution_SET_ONE} for 
    {\em LSUbl16dep, MlaFEDB, LetB1} and {\em TMEM16F}.
     }
    \label{Fig_viewing_angle_distribution_SET_TWO}
\end{figure}

\begin{figure}
    \centering
    \includegraphics[width=5.5in]{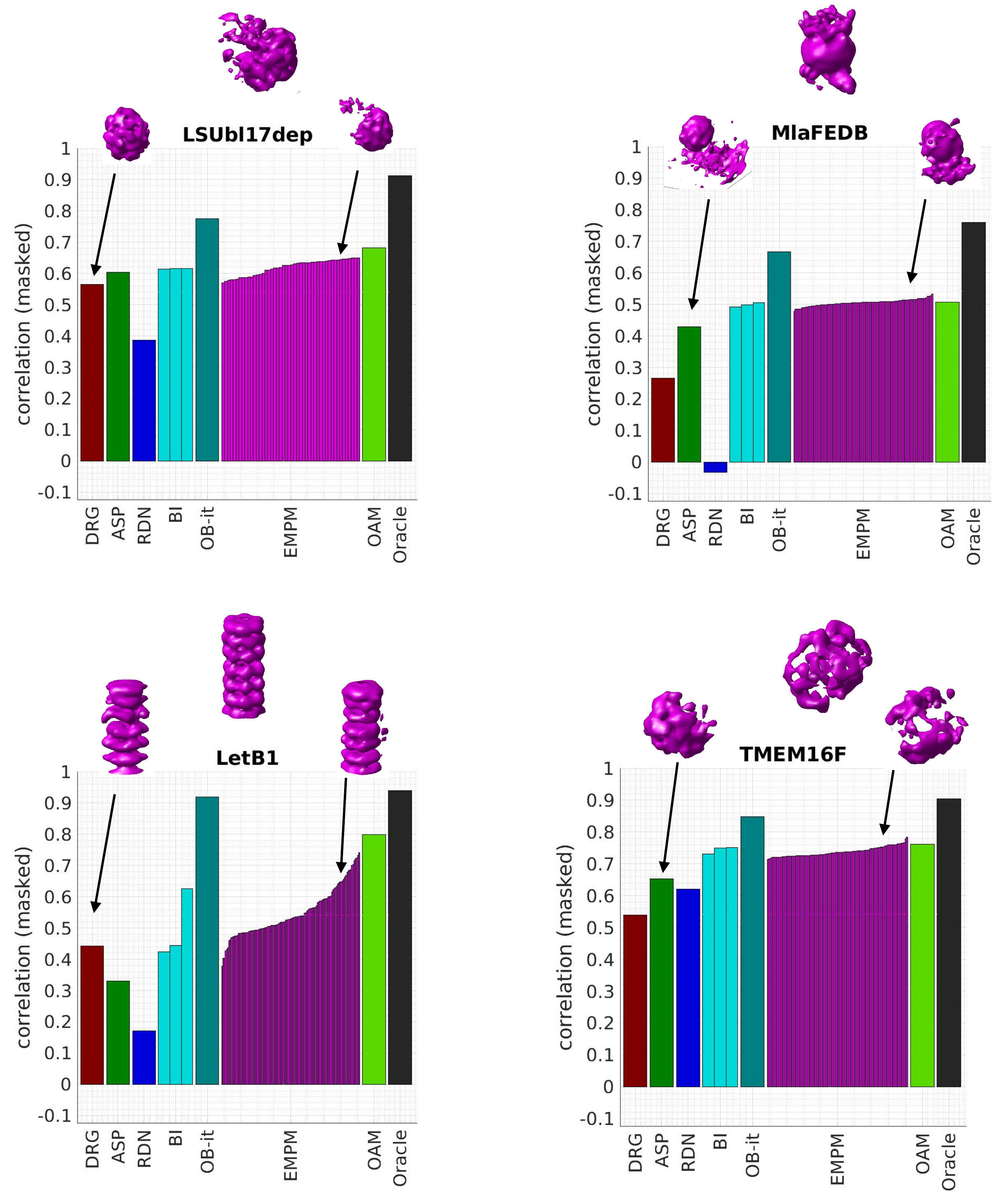}
    \caption{
{
    As in Fig. \ref{Fig_collect_SET_ONE}.
    we illustrate the aligned correlation $\cZ^{\masked}$ between our various reconstructions and the ground-truth. Here we examine reconstructions for
  {\em LSUbl16dep, MlaFEDB, LetB1} and {\em TMEM16F}.
    The different horizontal bars indicate the values of $\cZ^{\masked}$ for different reconstruction strategies (coded by color) with different trials from the same strategy adjacent to one another (sorted by $\cZ^{\masked}$).
    }
{At the top of each panel (centered) is the ``Oracle" low resolution
reconstruction for the indicated molecule. To the right in each panel, we show an EMPM 
reconstruction and to the left, we show a reconstruction obtained using one of the other pipelines
which achieves poorer correlation scores. Note that the $\cZ^{\masked}$ values
do tend to correlate with the image reconstruction quality.
{\em letB1} has a highly non-uniform distribution of ground-truth Euler angles. The
distributions in the other cases are more uniform. (See Fig. \ref{Fig_view_angle}.)
}}
    \label{Fig_collect_SET_TWO}
\end{figure}

\section{EMPM parameters:}
\label{sec_empm_param}

\begin{figure}
    \centering
    \includegraphics[width=6.5in]{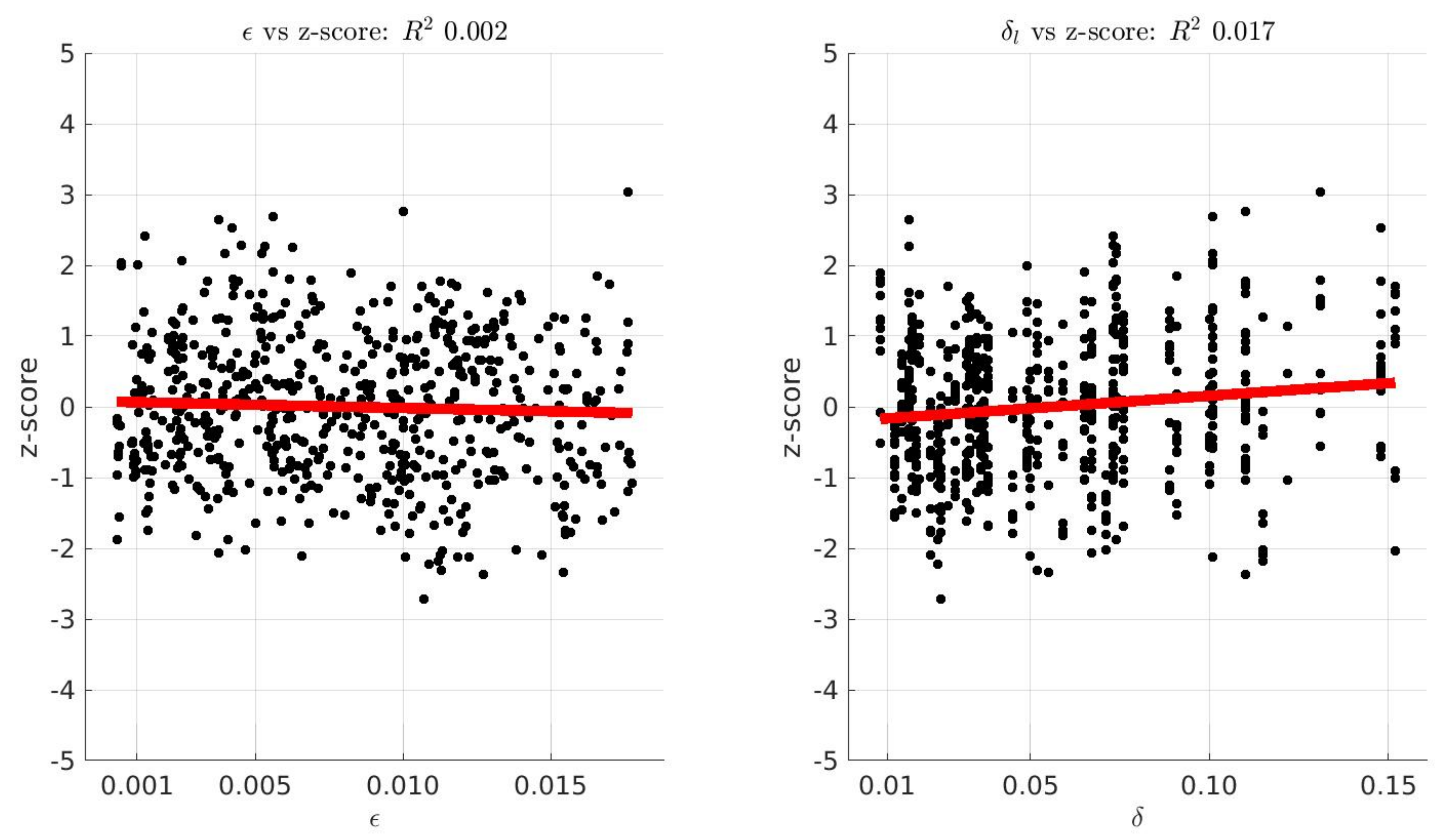}
    \caption{
{
    In this figure we illustrate the sensitivity of our results to perturbations in $\epsilon_{\tol}$ (left) and $\delta_{\local}$ (right).
    In each case we collect the results from all our reconstructions (across 8 molecules), transform each $\cZ^{\masked}$ into a normalized z-score (i.e., subtracting the mean $\cZ^{\masked}$ for that molecule, and then dividing by the standard-deviation across reconstructions for that molecule), and then plot the normalized z-score as a function of the selected parameter.
    Each dot corresponds to a reconstruction (633 in total), with the red line illustrating a linear fit.
    While there is a slight benefit to larger values of $\delta_{\local}$, and an even slighter benefit to smaller values of $\epsilon_{\tol}$, the variation in quality across reconstructions is more than an order of magnitude larger than the variation that can be attributed to the parameters (see $R^{2}$ values).
    }}
    \label{Fig_pm_collect_sensitivity_tolerance_sigma_FIGAB_strip}
\end{figure}

\begin{figure}
    \centering
    \includegraphics[width=3.5in]{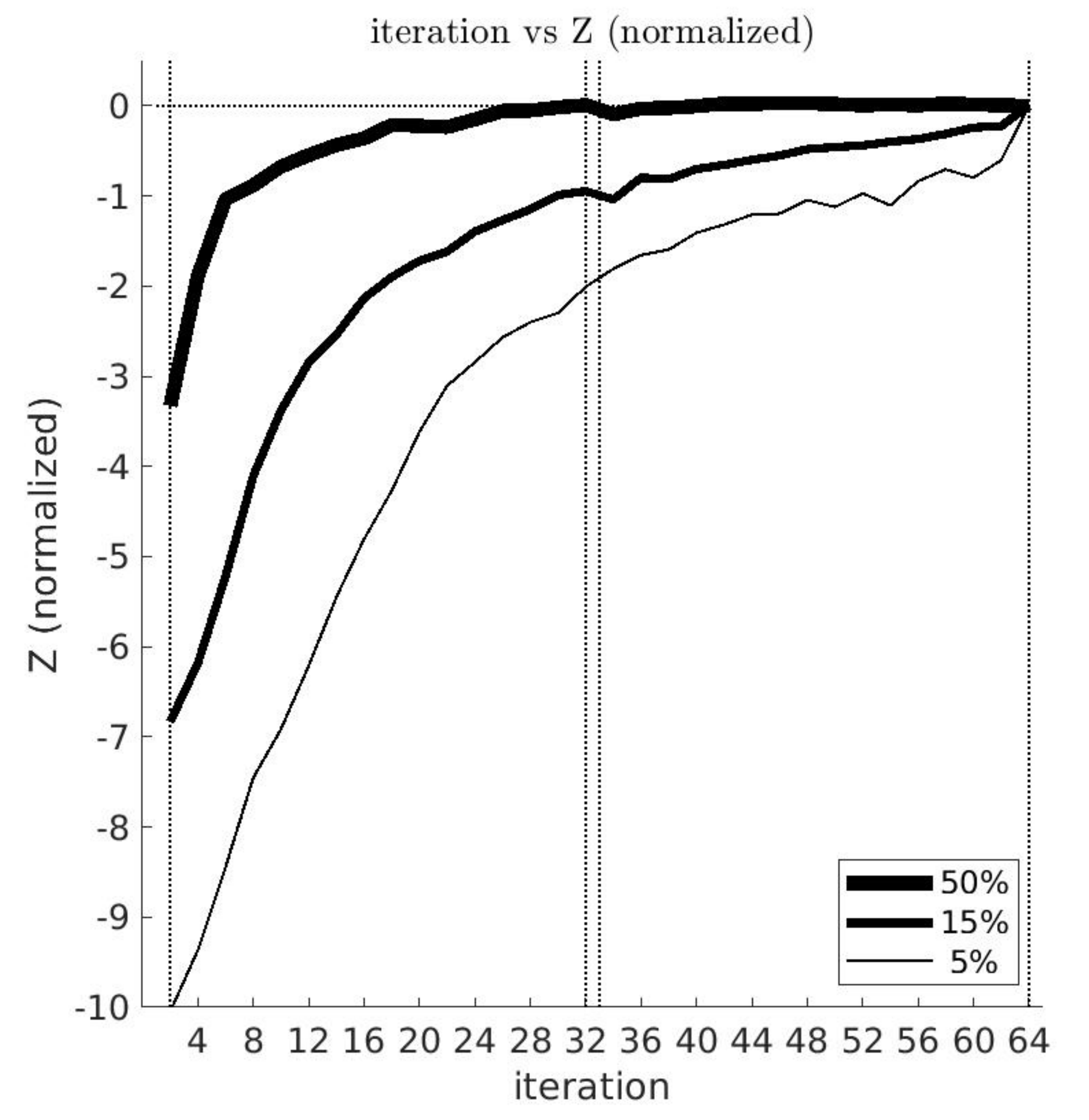}
    \caption{
{
    In this figure we illustrate the sensitivity of our results to the number of iterations chosen. For this figure we collect the results from all our reconstructions (across 8 molecules), then transform each (unmasked) $\cZ$ into a normalized score (i.e., subtracting the final (unmasked) $\cZ$ for that reconstruction, and then dividing by the standard-deviation of the final $\cZ$ taken across reconstructions for that molecule).
    We then plot the $50^{\th}$, $15^{\th}$ and $5^{\th}$ percentile for the distribution of normalized scores as a function of the even iterations (i.e., those corresponding to an entropy-maximizing alignment).
    Note that the former 32 iterations correspond to `phase-one', where the principal-mode compression uses the empirical principal-modes, while the latter 32 iterations correspond to `phase-two', where the principal-mode compression uses the volumetric principal-modes estimated after phase-one.
    Note that the median performance plateaus well before the final iteration of phase-two.
    If we were to halt our program earlier (after, say, only 16 iterations of phase-two) our results do not change considerably.
    }}
    \label{Fig_pm_collect_sensitivity_iteration_FIGC_strip}
\end{figure}

{
When performing the EMPM reconstructions shown in the main text, we fixed $\kmax=48$, corresponding to a maximum resolution of $\sim 20$\AA, as is typical for low-resolution estimates. Our EMPM algorithm also requires two other user-specified parameters: the overall tolerance $\epsilon_{\tol}$ and the width of the local displacement disk $\delta_{\local}$. As shown in Fig \ref{Fig_pm_collect_sensitivity_tolerance_sigma_FIGAB_strip}, the typical quality of our results does not change appreciably if these parameters are varied across an order of magnitude. We also run our algorithm for $64$ iterations. As shown in Fig \ref{Fig_pm_collect_sensitivity_iteration_FIGC_strip}, the median quality of our results does not deteriorate appreciably if we reduce the number of iterations by a modest
factor. }

{
As mentioned in the main text, the overall computation time for our algorithm is roughly constant per iteration (i.e., between $30-50$ minutes on a desktop workstation for all $64$ iterations). In principle, we expect the overall computation time of our algorithm to grow as $\epsilon_{\tol}$ becomes very small or $\delta_{\local}$ becomes very large. The reason we do not see a strong dependence on computation time for the range of parameters shown here is that, within this regime, the array dimensions involved are not so large (i.e., $\sim 10^{2}-10^{3}$ per dimension). As a consequence, roughly half of the overall computational effort is spent on memory movement rather than floating-point operations.
}

\begin{figure}
    \centering
    \includegraphics[width=5.5in]{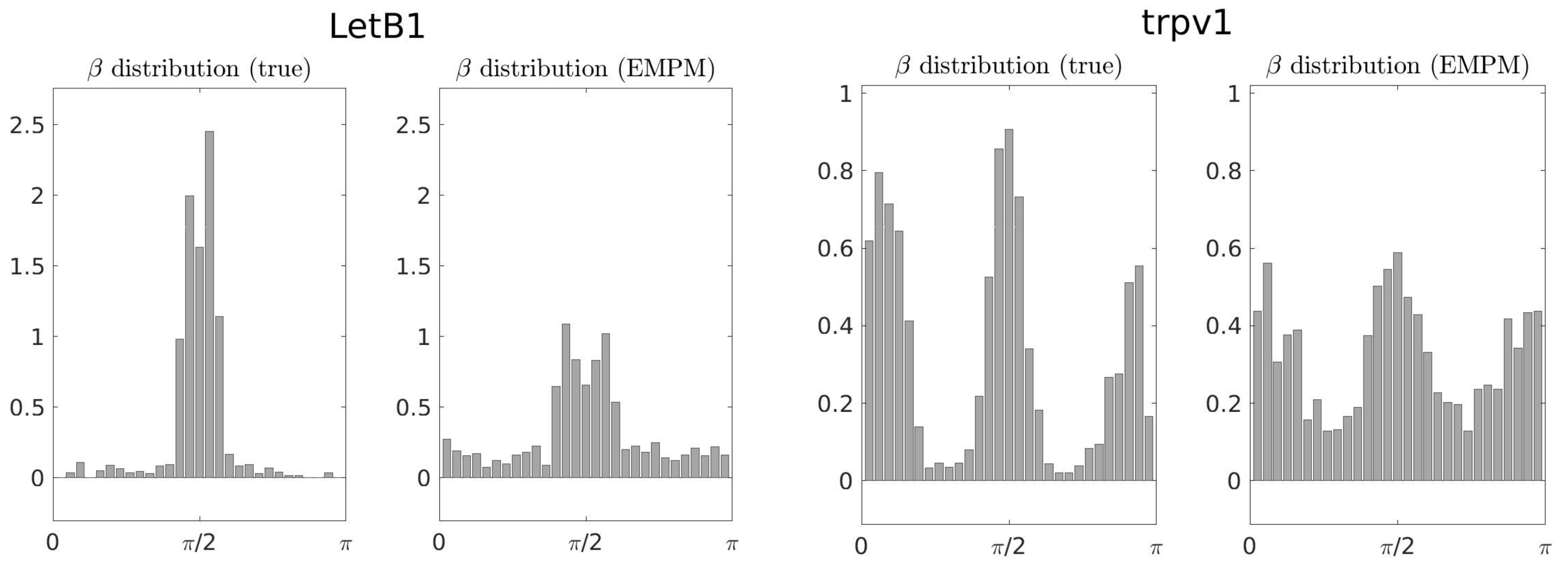}
    \caption{
{
    We illustrate the distribution of polar viewing-angle $\beta$ (ranging from $[0,\pi]$) for the images used in reconstructing {\em LetB1} (left) and {\em trpv1} (right). 
In each case, the left-hand histogram
shows the ground-truth distribution, while the right-hand histogram shows the distribution estimated from the reconstructions shown earlier.
This distribution is normalized by $\sin(\beta)$ so that a uniform distribution of viewing-angles on the sphere will appear constant as a function of $\beta$ (rather than being proportional to $\sin(\beta)$).
    Note that, even after this normalization, the viewing-angles for {\em LetB1} and
{\em trpv1} are strongly concentrated in certain regions of the sphere.
The ground-truth distribution of $\beta$ for {\em LetB1}, for example, 
has more than $90\%$ of its mass concentrated within $\pi/14$ of the equator.
    }}
    \label{Fig_view_angle}
\end{figure}


\begin{thebibliography}{10}

\bibitem{Milne2013}
Jacqueline L~S Milne, Mario~J. Borgnia, Alberto Bartesaghi, Erin E~H Tran,
  Lesley~A. Earl, David~M. Schauder, Jeffrey Lengyel, Jason Pierson, Ardan
  Patwardhan, and Sriram Subramaniam.
\newblock {Cryo-electron microscopy - A primer for the non-microscopist}.
\newblock {\em FEBS Journal}, 280:28--45, 2013.

\bibitem{Cheng2015}
Yifan Cheng, Nikolaus Grigorieff, Pawel~A. Penczek, and Thomas Walz.
\newblock A primer to single-particle cryo-electron microscopy.
\newblock {\em Cell}, 161:439--449, 2015.

\bibitem{Epstein}
C.~L. Epstein.
\newblock {\em Introduction to the Mathematics of Medical Imaging}.
\newblock SIAM, 2008.

\bibitem{Natterer}
F.~Natterer.
\newblock {\em The Mathematics of Computerized Tomography}.
\newblock SIAM, 2001.

\bibitem{Kimanius2016}
Dari Kimanius, Bjorn~O Forsberg, Sjors~HW Scheres, and Erik Lindehl.
\newblock {Accelerated cryo-EM structure determination with parallelisation
  using GPUs in RELION-2}.
\newblock {\em eLife}, 5:e18722, 2016.

\bibitem{Kimaniusetal2021}
D.~Kimanius, L.~Dong, G.~Sharov, T.~Nakane, and Scheres S.H.W.
\newblock New tools for automated cryo-em single-particle analysis in
  relion-4.0.
\newblock {\em Biochem}, 478:4169--4185, 2021.

\bibitem{Scheres2012}
Sjors H~W Scheres.
\newblock A {B}ayesian view on {cryo-EM} structure determination.
\newblock {\em J. Mol. Biol.}, 415:406--418, 2012.

\bibitem{ZBBD21}
E.D. Zhong, T.~Bepler, B.~Berger, and J.H. Davis.
\newblock Cryodrgn: reconstruction of heterogeneous cryo-em structures using
  neural networks.
\newblock {\em Nature Methods}, 18:176--185, 2021.

\bibitem{AspireWebsite}
Aspire software, 2022.
\newblock {\tt http://spr.math.princeton.edu}.

\bibitem{Punjani2017}
A.~Punjani, J.L. Rubinstein, D.J. Fleet, and M.A. Brubaker.
\newblock {cryoSPARC}: algorithms for rapid unsupervised {cryo-EM} structure
  determination.
\newblock {\em Nat. Methods}, 14:290--296, 2017.

\bibitem{Scheres2012b}
Sjors H~W Scheres.
\newblock {RELION: Implementation of a Bayesian approach to cryo-EM structure
  determination}.
\newblock {\em J. Struct. Biol.}, 180(3):519--530, 2012.

\bibitem{DTCPLW16}
J.H. Davis, Y.Z. Tan, B.~Carragher, C.S. Potter, D.~Lyumkis, and J.R.
  Williamson.
\newblock Modular assembly of the bacterial large ribosomal subunit.
\newblock {\em Cell}, 167:1610--1622, 2016.

\bibitem{rangancode}
{E}{M}{P}{M} software, 2023.
\newblock {\tt https://github.com/adirangan/dir\_cryoem}.

\bibitem{bracewell}
R~Bracewell.
\newblock {\em The {F}ourier Transform and Its Applications}.
\newblock McGraw-Hill, 3rd edition, 1999.

\bibitem{Zhao2014}
Zh. Zhao and A.~Singer.
\newblock Rotationally invariant image representation for viewing direction
  classification in {cryo-EM}.
\newblock {\em J. Struct. Biol.}, 186:153--166, 2014.

\bibitem{Barnett2017}
A.~Barnett, L.~Greengard, A.~Pataki, and M.~Spivak.
\newblock Rapid solution of the {cryo-EM} reconstruction problem by frequency
  marching.
\newblock {\em SIAM J. Imaging Sci.}, 10(3):1170--1195, 2017.

\bibitem{RSAB20}
Aaditya Rangan, Marina Spivak, Joakim And{\'e}n, and Alex Barnett.
\newblock Factorization of the translation kernel for fast rigid image
  alignment.
\newblock {\em Inverse Problems}, 36(2), 2020.

\bibitem{Zhao2016}
Zh. Zhao, Y.~Shkolnisky, and A.~Singer.
\newblock Fast steerable principal component analysis.
\newblock {\em IEEE Trans. Comput. Imaging}, 2(1):1--12, 2016.

\bibitem{finufft}
A~H Barnett, J~F Magland, and L~af~Klinteberg.
\newblock A parallel non-uniform fast {F}ourier transform library based on an
  ``exponential of semicircle'' kernel, 2019.

\bibitem{Sigworth1998}
F~J Sigworth.
\newblock {A maximum-likelihood approach to single-particle image refinement}.
\newblock {\em J. Struct. Biol.}, 122(3):328--39, 1998.

\bibitem{Scheres2009}
S~H~W Scheres, M~Valle, P~Grob, E~Nogales, and J.-M. Carazo.
\newblock Maximum likelihood refinement of electron microscopy data with
  normalization errors.
\newblock {\em J. Struct. Biol.}, 166(2):234--240, 2009.

\bibitem{Sigworth2010}
F.~J. Sigworth, Doerschuk P.C., J.-M. Carazo, and S.H.W. Scheres.
\newblock {An introduction to maximum-likelihood methods in Cryo-EM}.
\newblock In {\em Methods in Enzymology. Cryo-EM, Part B: 3D reconstruction},
  pages 263--294. Academic Press., 2010.

\bibitem{Lyumkis2013}
Dmitry Lyumkis, Axel~F. Brilot, Douglas~L. Theobald, and Nikolaus Grigorieff.
\newblock {Likelihood-based classification of cryo-EM images using FREALIGN}.
\newblock {\em J. Struct. Biol.}, 183(3):377--388, 2013.

\bibitem{Baxter2009}
W.~Baxter, R.~Grassucci, Haixiao Gao, and J.~Frank.
\newblock Determination of signal-to-noise ratios and spectral snrs in cryo-em
  low-dose imaging of molecules.
\newblock {\em Journal of structural biology}, 166 2:126--32, 2009.

\bibitem{Rangan22}
Aaditya~V Rangan.
\newblock Radial-recombination for rigid rotational alignment of images and
  volumes.
\newblock {\em Inverse Problems (submitted)}, 2022.

\bibitem{bendory_mra_2021}
Tamir Bendory, Ariel Jaffe, William Leeb, Nir Sharon, and Amit Singer.
\newblock {Super-resolution multi-reference alignment}.
\newblock {\em Information and Inference: A Journal of the IMA},
  11(2):533--555, 02 2021.

\bibitem{bendory_mra_2022}
Noam Janco and Tamir Bendory.
\newblock An accelerated expectation-maximization algorithm for multi-reference
  alignment.
\newblock {\em IEEE Transactions on Signal Processing}, 70:3237--3248, 2022.

\bibitem{perry_mra_2019}
Amelia Perry, Jonathan Weed, Afonso~S. Bandeira, Philippe Rigollet, and Amit
  Singer.
\newblock The sample complexity of multireference alignment.
\newblock {\em SIAM Journal on Mathematics of Data Science}, 1(3):497--517,
  2019.

\end{thebibliography}
\end{document}